\documentclass[11pt]{article}
\usepackage{amscd}
\usepackage{amsfonts}
\usepackage{amsmath}
\usepackage{amssymb}
\usepackage{amsthm}
\usepackage{bbm}
\usepackage{CJK}
\usepackage{fancyhdr}
\usepackage{graphicx}
\usepackage{hyperref}
\usepackage{indentfirst}
\usepackage{latexsym}
\usepackage{mathrsfs}
\usepackage{xypic}

\usepackage[top=1in,bottom=1in,left=1.25in,right=1.25in]{geometry}
\def\<{\langle}
\def\>{\rangle}
\def\a{\alpha}
\def\b{\beta}

\def\c{\cdot}

\def\D{\Delta}

\def\e{\eta}
\def\g{\gamma}

\def\o{\otimes}

\def\ra{\rightarrow}

\def\v{\epsilon }

\date{}
\begin{document}
\renewcommand{\baselinestretch}{1.2}
\renewcommand{\arraystretch}{1.0}
\title{\bf Central invariants and enveloping algebras \\of braided Hom-Lie algebras}
\author{{\bf Shengxiang Wang$^{1}$,Xiaohui Zhang$^{2}$
        and Shuangjian Guo$^{3}$\footnote
        {Corresponding author(Shuangjian Guo): shuangjianguo@126.com} }\\
1.~ School of Mathematics and Finance, Chuzhou University,\\
 Chuzhou 239000,  China \\
2.~  School of Mathematical Sciences, Qufu Normal University, \\Qufu Shandong 273165, China.\\
 3.~ School of Mathematics and Statistics, Guizhou University of\\ Finance and Economics,  Guiyang 550025, China.}
 \maketitle
\begin{center}
\begin{minipage}{13.cm}

{\bf \begin{center} ABSTRACT \end{center}}
Let $(H,\alpha)$ be a monoidal Hom-Hopf algebra
and $^{H}_{H}\mathcal{HYD}$ the Hom-Yetter-Drinfeld category over $(H,\alpha)$.
Then in this paper, we first introduce the definition of braided Hom-Lie algebras
and show that each monoidal Hom-algebra in  $^{H}_{H}\mathcal{HYD}$ gives rise to a braided Hom-Lie algebra.
Second, we prove that if $(A,\beta)$ is a sum of two $H$-commutative monoidal Hom-subalgebras,
then the commutator Hom-ideal $[A,A]$ of $A$ is nilpotent.
Also, we study the central invariant of braided Hom-Lie algebras as a generalization of generalized Lie algebras.
Finally, we obtain a construction of the enveloping algebras of braided Hom-Lie algebras and show that the enveloping algebras
are $H$-cocommutative Hom-Hopf algerbas.\\


 \end{minipage}
 \end{center}
 \normalsize\vskip1cm

\section*{Introduction}
\def\theequation{0. \arabic{equation}}
\setcounter{equation} {0}

Hom-algebras were first introduced in the Lie algebra setting \cite{Hartwig} with motivation from physics
though its origin can be traced back in earlier literature such as \cite{Hu}.
In a Hom-Lie algebra, the Jacobi identity is replaced by the so called Hom-Jacobi identity via a homomorphism.
In 2008, Makhlouf and Silvestrov \cite{Makhlouf2008} introduced the definition of Hom-associative algebras,
where the associativity of a Hom-algebra is twisted by an endomorphism (here we call it the Hom-structure map).
The generalized notions, including Hom-bialgebras, Hom-Hopf algebras were developed
in \cite{Dekkar}, \cite{Makhlouf2009}, \cite{Makhlouf2010}, \cite{AS3}.
Further research on  Hom-Hopf algebras could be found in \cite{CWZ2014}, \cite{Gohr2010}, \cite{Makhlouf2018}, \cite{Yau2009},
\cite{zz} and references cited therein.

In \cite{Caenepeel2011}, Caenepeel and Goyvaerts studied Hom-Lie algebras and Hom-Hopf algebras from a categorical view point,
 they proved a (co)monoid in the Hom-category is a Hom-(co)algebra, and a bimonoid in the Hom-category is a monoidal Hom-bialgebra.
  Note that a monoidal Hom-Hopf algebra is a Hom-Hopf algebra if and only if the Hom-structure map is involutional.
Later, Graziani et al.  \cite{Graziani}  defined BiHom-Hopf algebras using two commuting multiplicative linear maps $\alpha,\beta$,
unified Hom-Hopf algebras and monoidal Hom-Hopf algebras by setting $\alpha=\beta$ and $\alpha=\beta^{-1}$ respectively.

Recently, the theory of Hom-Yetter-Drinfeld categories have attracted attention in mathematics and mathematical physics.
In \cite{Makhlouf2014}, Makhlouf and Panaite defined Yetter-Drinfeld modules over Hom-bialgebras
 and shown that Yetter-Drinfeld modules over a Hom-bialgebra with bijective structure map provide solutions of the Hom-Yang-Baxter equation.
Also Liu and Shen \cite{LIU2014}, Chen and Zhang  \cite{chen2014} studied Hom-Yetter-Drinfeld modules over monoidal Hom-bialgebras in a slightly different way to \cite{Makhlouf2014}.
As a part of the theory of Hom-Yetter-Drinfeld categories,
we \cite{wang&GUO2014} gave sufficient and necessary conditions for the Hom-Yetter-Drinfeld category $^{H}_{H}\mathcal{HYD}$
to be symmetric and pseudosymmetric respectively.
With the symmetries of Hom-Yetter-Drinfeld categories, it is a natural question to ask whether we can extend the notion of
monoidal Hom-Lie algebras to  Hom-Yetter-Drinfeld categories.
This becomes our first motivation of writing this paper.

It is well known that  Lie algebras in braided monoidal categories is a very important part of Lie theories.
As a generalization of Lie superalgebras \cite{Kac}  and Lie color algebras \cite{Scheunert},
 Manin \cite{Manin} studied Lie algebras in some symmetric categories from an algebraic point of view.
Later, Cohen, Fishman and Westreich \cite{Cohen1994} studied Lie algebras in the category of modules over triangular Hopf algebras
and proved Schur's double centralizer theorem,
 Fishman and Montgomery \cite{Fishman1994} did similar work in the category of comodules over cotriangular Hopf algebras.
Later, Bahturin, Fishman and Montgomery \cite{Bahturin} studied the structure of the
 generalized Lie algebras in the category of comodules.

Wang \cite{Wang2002} introduced the notion of  generalized Lie algebras in Yetter-Drinfeld categories and
extended the Kegel's theorem to generalized Lie algebras.
Later, we \cite{Wang2014} extended Wang's results in \cite{Wang2002} to Hom-Lie algebras in Yetter-Drinfeld categories,
which unifies the notions of Hom-Lie superalgebras in \cite{Ammar} and Hom-Lie color algebras in \cite{Yuan}.
In the present paper, we will study monoidal Hom-Lie algebras in Hom-Yetter-Drinfeld categories, which is different from \cite{Wang2014} in two aspects.
First, Hom-Yetter-Drinfeld categories include Yetter-Drinfeld categories as a special case.
Second, the main purpose of this paper is to study the central invariants an enveloping algebras of  braided Hom-Lie algebras,
 which has not been involved in \cite{Wang2014}.

 This paper is organized as follows.
In Section 1, we recall some basic definitions about monoidal Hom-Hopf algebras and Hom-Yetter-Drinfeld modules.

In Section 2, we define braided Hom-Lie algebras
and show that any monoidal Hom-algebra in  $^{H}_{H}\mathcal{HYD}$ gives rise to a braided
 Hom-Lie algebra by the natural bracket product (see Proposition 2.2),
and prove that if $(A,\beta)$ is H-semisimple and a sum of two $H$-commutative monoidal Hom-subalgebras,
 then $(A,\beta)$ is H-commutative (see Corollary 2.9).
In Section 3, we consider the central invariant of braided Hom-Lie algebras
(see Theorem 3.7).
In Section 4, we construct the  enveloping algebras of  braided Hom-Lie algebras and present its Hopf structures.
As an application, we study the  enveloping algebras of $End(V)$ and construct a Radford's Hom-biproduct
$(U(End(V))_{\sharp}^{\times} H, \delta\o id)$ (see Proposition 4.10).

\section{Preliminaries}
\def\theequation{\arabic{section}.\arabic{equation}}
\setcounter{equation} {0}

In this section, we recall some basic definitions and results related to our paper.
Throughout the paper, all algebraic systems are supposed to be over a field ${k}$.
The reader is referred to Caenepeel and Goyvaerts \cite{Caenepeel2011}
as general references about monoidal Hom-algebras and monoidal Hom-Lie algebras,
to Sweedler \cite{Sweedler} about Hopf algebras and Liu and Shen \cite{LIU2014} about Hom-Yetter-Drinfeld categories.

If $C$ is a coalgebra, we use the Sweedler-type notation for the comultiplication: $\Delta(c)=c_{1}\o c_{2}$,
for all $c\in C,$
 in which we often omit the summation symbols for convenience.
\medskip

{\bf 1.1 Hom-category}

Let $\mathcal{C}$ be a category. We introduce a new category
$\mathscr{H}(\mathcal{C})$ as follows: the objects are
couples $(X, \a_X)$, with $M \in \mathcal{C}$ and $\a_X \in
Aut_{\mathcal{C}}(X)$. A morphism $f: (X, \a_X)\rightarrow (Y, \a_Y)$
is a morphism $f : X\rightarrow Y$ in $\mathcal{C}$ such that $\a_Y
\circ f= f \circ \a_X$.

Specially, let $\mathscr{M}_k$ denote the category of $k$-spaces.
~$\mathscr{H}(\mathscr{M}_k)$ will be called the Hom-category
associated to $\mathscr{M}_k$. If $(X,\a_X) \in \mathscr{M}_k$, then
$\a_X: X\rightarrow X$ is obviously an isomorphism in
~$\mathscr{H}(\mathscr{M}_k)$. It is easy to show that
~$\widetilde{\mathscr{H}}(\mathscr{M}_k)$ =
(~$\mathscr{H}(\mathscr{M}_k),~\otimes,~(k, id),~\widetilde{a},
~\widetilde{l},~\widetilde{r}))$ is a monoidal category by
Proposition 1.1 in \cite{Caenepeel2011}:

$\bullet$ the tensor product of $(X,\a_X)$ and $(Y,
\a_Y)$ in ~$\widetilde{\mathscr{H}}(\mathscr{M}_k)$ is given by the
formula $(X, \a_X)\otimes (Y, \a_Y) = (X\otimes Y, \a_X \otimes \a_Y)$;

$\bullet$ for any $x \in X$, $y \in Y$, $z \in Z$, the associativity is given by the formulas
$$
\widetilde{a}_{X,Y,Z}((x\o y)\o z)=\a_X(x)\o (y\o \a_Z^{-1}(z));$$

$\bullet$ for any $x \in X$, $\lambda \in k$, the unit
constraints are given by the formulas
$$
\widetilde{l}_{X}(\lambda\o x)=\widetilde{r}_{X}(x\o \lambda)=\lambda\a_X(x).
$$

{\bf 1.2 Monoidal Hom-Hopf algebras}
\medskip

{\it Definition 1.1.} A {\it monoidal Hom-algebra} is an object $(A, \a)$ in the Hom-category
$\widetilde{\mathcal{H}}(\mathcal{M}_k)$ together with an element $1_A\in A$
and a linear map $m: A\o A\ra A,\,\, a\o b\mapsto ab$
such that
 \begin{eqnarray}
&&\a(a)(bc)=(ab)\a(c),~\a(ab)=\a(a)\a(b),\\
&&a1_A=1_A a= \a(a),~ \a(1_A)=1_A,
\end{eqnarray}
for all $a, b, c\in A$.
\smallskip

As noted in \cite{Caenepeel2011}, the definition of monoidal Hom-algebras is different from the definition of
Hom-associative algebras defined in \cite{Makhlouf2010}.
Specifically, the unitality condition in \cite{Makhlouf2010} is the usual untwisted one: $a1_A=1_A a= a$, for any $a\in A$,
and the condition (1.2) is not desired there.
\smallskip

{\it Definition 1.2.} A {\it monoidal Hom-coalgebra} is an object $(C, \g)$ in the category
$\widetilde{\mathcal{H}}(\mathcal{M}_k)$ together with linear maps
$\D: C\ra C\o C,\, \D(c)=c_1\o c_2$ and $\v: C \ra k$ such that
 \begin{eqnarray}
&&\g^{-1}(c_1)\o \D(c_2)=\D(c_1)\o \g^{-1}(c_2),~\D(\g(c))=\g(c_1)\o \g(c_2),\\
&&c_1\v(c_2)=\v(c_1)c_2=\g^{-1}(c), ~ \v(\g(c))=\v(c),
\end{eqnarray}
for all $c\in C$.
\smallskip

The definition of monoidal Hom-coalgebras is different from the definition of
Hom-coassociative coalgebras defined in \cite{Makhlouf2010}. The coassociativity condition is twisted by some endomorphism,
not necessarily by the inverse of the automorphism $\gamma$.
The counitality condition in \cite{Makhlouf2010} is the usual untwisted one: $c_1\v(c_2)=\v(c_1)c_2=c$, for any $c\in C$,
and the condition (1.5) is not needed there.
\medskip

{\it Definition 1.3.} A {\it monoidal Hom-bialgebra} $H=(H, \a, m,
1_H, \D, \v)$ is a bialgebra in the category
$\widetilde{\mathcal{H}}(\mathcal{M}_k)$. This means that $(H, \a,
m, 1_H)$ is a monoidal Hom-algebra and $(H, \a, \D, \v)$ is a
monoidal Hom-coalgebra such that $\D$ and $\v$ are Hom-algebra maps,
that is, for any $h, g\in H$,
 \begin{eqnarray*}
\D(hg)=\D(h)\D(g),& ~ \D(1_H)=1_H\o 1_H,\\
\v(hg)=\v(h)\v(g),& ~ \v(1_H)=1_k.
\end{eqnarray*}

A monoidal Hom-bialgebra $(H, \a)$ is called a {\it monoidal
Hom-Hopf algebra} if there exists a morphism (called the antipode) $S:
H\ra H$ in $\widetilde{\mathcal{H}}(\mathcal{M}_k)$ (i.e. $S\circ
\a=\a \circ S$), which is the convolution inverse of the
identity morphism $id_H$ (i.e. $S* id_H =\eta_H\circ \v_H=id_H *
S$), this means for any $h\in H$,
 \begin{eqnarray}
S(h_1)h_2=\v(h)1_H=h_1S(h_2).
\end{eqnarray}


{\bf 1.3  Hom-Yetter-Drinfeld categories}
\medskip

{\it Definition 1.4.} Let $(A, \a)$ be a monoidal Hom-algebra. A {\it left $(A,
 \a)$-Hom-module}
  consists of $(M, \mu)\in \widetilde{\mathcal{H}}(\mathcal{M}_k)$
  together with a morphism $\psi: A\o M \ra M,\, \psi(a\o m)=a\cdot m$ such
  that
 \begin{eqnarray}
 \a(a)\cdot(b\cdot m)=(ab)\cdot\mu(m), ~ 1_A \cdot m=\mu(m),
 \mu(a\cdot m)=\a(a)\cdot\mu(m),
\end{eqnarray}
 for all $a, b \in A$ and $m\in M$.
  \smallskip

A morphism $f: M\ra N$ is called {\it left
 $A$-linear} if
 $f(am)=a f(m)$, for any $a\in A, m\in M$ and $f\circ \mu=\nu\circ
 f$.

{\it Definition 1.5.} Let $(C, \g)$ be a monoidal Hom-coalgebra. A
{\it left $(C,
 \g)$-Hom-comodule}
  consists of $(M, \mu)\in \widetilde{\mathcal{H}}(\mathcal{M}_k)$
  together with a morphism $\rho_M: M\ra C\o M,\, \rho_M(m)=m_{(-1)}\o m_{0}$ such
  that
 \begin{eqnarray}
 &&\D_{C}(m_{(-1)})\o \mu^{-1}(m_{0})=\g^{-1}(m_{(-1)})\o (m_{0(-1)}\o m_{00}),\\
  && \rho_M(\mu(m))=\g(m_{(-1)})\o \mu(m_{0}),\,\v(m_{(-1)})m_{0}=\mu^{-1}(m),
\end{eqnarray}
 for all $m\in M$.
  \smallskip

Let $(M, \mu)$ and $(N,
 \nu)$ be two left $(C,
 \g)$-Hom-comodules. A morphism $g: M\ra N$ is called {\it left
 $C$-colinear} if $g\circ \mu=\nu\circ
 g$ and
 $m_{(-1)}\o g(m_{0})=g(m)_{(-1)}\o g(m)_{0}$, for any $m\in M$.

{\it Definition 1.6.} Let $(H, \alpha)$ be a monoidal Hom-Hopf algebra. A
{\it left-left $(H, \alpha)$-Hom-Yetter-Drinfeld module}
 is an object $(M, \beta)\in \widetilde{\mathcal{H}}(\mathcal{M}_k)$,
 such that $(M, \beta)$ is both a left $(H, \alpha)$-Hom-module
 and a left $(H, \alpha)$-Hom-comodule with the following compatibility condition:
 \begin{eqnarray}
\rho(h \cdot m)=(h_{11}\alpha^{-1}(m_{(-1)}))S(h_{2})\otimes \alpha(h_{12})\cdot m_{0},
\end{eqnarray}
 for all $h\in H$ and $m\in M$.

One has that Eq. (1.9) is equivalent to the following equation:
$$
h_1 m_{(-1)}\otimes h_2\cdot m_0=(h_1\cdot \beta^{-1}(m))_{(-1)}h_2
\otimes \beta((h_1\cdot \beta^{-1}(m))_{0}).
$$

{\it Definition 1.7.} Let $(H, \alpha)$ be a monoidal Hom-Hopf algebra. A
{\it Hom-Yetter-Drinfeld category} $^{H}_{H}\mathcal{HYD}$
 is a braided monoidal category whose objects are left-left $(H, \alpha)$-Hom-Yetter-Drinfeld modules, morphisms are both
 left $(H, \alpha)$-linear and $(H, \alpha)$-colinear maps, and its braiding $C_{-, -}$  is given by
\begin{eqnarray}
 C_{M, N} (m \o n) = m_{(-1)}\c \nu^{-1}(n)\o \mu(m_{(0)}),
\end{eqnarray}
for all $m \in (M,\mu)\in{}^{H}_{H}\mathcal{HYD}$ and $n \in (N,\nu) \in{}^{H}_{H}\mathcal{HYD}$.
\medskip

{\it Definition 1.8.} Let $(A,\beta)$ be an object in $^{H}_{H}\mathcal{HYD}$,
 the braiding $C$ is called symmetric on $A$ if the following condition holds:
 \begin{eqnarray}
a_{(-1)}\cdot\beta^{-1}(b)\otimes\beta(a_{0})=\beta(b_{0})\otimes S^{-1}(b_{(-1)})\cdot\beta^{-1}(a);
\end{eqnarray}
$A$ is called $H$-commutative if
 \begin{eqnarray}
(a_{(-1)}\cdot\beta^{-1}(b))\beta(a_{0})=ab,
\end{eqnarray}
$A$ is called $H$-cocommutative if
 \begin{eqnarray}
a_{1(-1)}\cdot\beta^{-1}(a_{2})\o\beta(a_{10})=a_{2}\o a_{1},
\end{eqnarray}
for all $a,b\in A.$

\section{Braided Hom-Lie algebras}
\def\theequation{\arabic{section}. \arabic{equation}}
\setcounter{equation} {0}

In this section, we first introduce the concept of braided Hom-Lie algebras and
 show that each monoidal Hom-algebra in $^{H}_{H}\mathcal{HYD}$ gives rise to a  braided Hom-Lie algebras.
Also we study the braided Lie structures of monoidal Hom-algebras in $^{H}_{H}\mathcal{HYD}$
as a generalization of results in \cite{Bahturin}, \cite{Wang2002} and \cite{Wang2014}.
\medskip

From now on, we always assume that $(H,\alpha)$ is a monoidal Hom-Hopf algebra
 and  $^{H}_{H}\mathcal{HYD}$ the  Hom-Yetter-Drinfeld category over  $(H,\alpha)$.
 \medskip

\noindent{\bf Definition 2.1.}
A monoidal Hom-Lie algebra in $^{H}_{H}\mathcal{HYD}$,
called by a braided Hom-Lie algebra,
 is a triple $(L,[,],\beta)$,
where $L$ is an object in $^{H}_{H}\mathcal{HYD}$,
$\beta:L\rightarrow L$ is a homomorphism in $^{H}_{H}\mathcal{HYD}$
and $[,]:L\o L\rightarrow L$ is a morphism in $^{H}_{H}\mathcal{HYD}$ satisfying

(i)~Braided Hom-skew-symmetry:
\begin{eqnarray}
[l,l']=-[l_{(-1)}\c \beta^{-1}(l'),\beta(l_{0})],~l,l'\in L.
\end{eqnarray}

(ii)~Braided Hom-Jacobi identity:
\begin{eqnarray}
\{l\o l'\o l''\}+\{(C\o 1)(1\o C)(l\o l'\o l'')\}+\{(1\o C)(C\o 1)(l\o l'\o l'')\}=0,
\end{eqnarray}
for all $l,l',l''\in L$, where $\{l\o l'\o l''\}$ denotes $[\beta(l),[l',l'']]$ and $C$ the braiding for $L$.
\medskip

\noindent{\bf Proposition 2.2.}
Let $(A,\beta)$ be a monoidal Hom-algebra in $^{H}_{H}\mathcal{HYD}$.
Assume that the braiding $C$ is symmetric on $A$.
Then the triple $(A,[,],\beta)$ is a braided Hom-Lie algebra,
where the bracket product is defined by
\begin{eqnarray}
[,]:A\otimes A\rightarrow A~by~[a,b]=ab-(a_{(-1)}\cdot \beta^{-1}(b))\beta(a_{0}),
\end{eqnarray}
for all $a,b\in A$.

{\bf Proof.} Denote $A^{-}=(A,[,],\beta)$.
It is clear that the bracket product is a morphism in $^{H}_{H}\mathcal{HYD}$,
so it remains to verify that the conditions (i) and (ii) of Definition 2.1 hold.

For braided Hom-skew-symmetry, we have
$
[a_{(-1)}\cdot \beta^{-1}(b), \beta(a_{0})]
=(a_{(-1)}\cdot \beta^{-1}(b))\beta(a_{0})-((a_{(-1)}\cdot \beta^{-1}(b))_{(-1)}\cdot \beta(\beta^{-1}(a_{0})))\beta((a_{(-1)}\cdot \beta^{-1}(b))_{0})
=(a_{(-1)}\cdot \beta^{-1}(b))\beta(a_{0})-ab
=-[a,b],
$
as desired.
The last equality holds since the braiding $C$ is symmetric on $A$.

Similarly, one may check the braided Hom-Jacobi identity by the Hom-associativity of $A$ routinely.
And this finishes the proof. $\hfill \Box$
\medskip

\noindent{\bf  Example 2.3.}
 Let $(H, \alpha)$ be a commutative monoidal Hom-Hopf algebra.
By Example 4.3 in \cite{LIU2014}, $(H, \alpha)$  is a Hom-Yetter-Drinfeld module with left $(H, \alpha)$-action
$h\cdot g=(h_{1}\alpha^{-1}(g))S(\alpha(h_{2}))$ and left $(H, \alpha)$-coaction by the
Hom-comultiplication $\Delta$, note it by $H_{1}=(H_{1}, \mbox{adjoint},\Delta,\alpha)$.
By Corollary 5.4 in \cite{wang&GUO2014}, the braiding $C$ is symmetric on $H_{1}$,
 then $H_{1}^{-}$ is a braided Hom-Lie algebra.
\medskip

\noindent{\bf  Example 2.4.}
 Let $(H, \alpha)$ be a cocommutative monoidal Hom-Hopf algebra.
By Example 2.7 in \cite{wang&GUO2014}, $(H, \alpha)$  is a Hom-Yetter-Drinfeld module with left $(H, \alpha)$-action by the
 Hom-multiplication $m$ and left $(H, \alpha)$-coaction $\rho(h)=h_{11}\alpha^{-1}(S(h_{2}))\otimes \alpha(h_{12})$,
and note it by $H_{2}=(H_{2}, m,  \mbox{coadjoint},\alpha)$.
By Corollary 4.4 in \cite{wang&GUO2014}, the braiding $C$ is symmetric on $H_{2}$,
 then $H_{2}^{-}$ is a braided Hom-Lie algebra.
\medskip

\noindent{\bf  Example 2.5.}
Let $H=k\{1_{H},h\}$ be a monoidal Hom-Hopf algebra with an automorphism
 $\alpha:H\rightarrow H,\alpha(1_{H})=1_{H},\alpha(h)=-h, $
 where the Hom-algebra structure  is defined by
\begin{eqnarray*}
 1_{H}1_{H}=1_{H},1_{H}h=h1_{H}=-h,h^{2}=0,
\end{eqnarray*}
 the Hom-coalgebra structure  is defined by
\begin{eqnarray*}
 \Delta(1_{H})=1_{H}\o1_{H},
 \Delta(h)=(-h)\o1_{H}+1_{H}\o(-h), \epsilon(1_{H})=1,\epsilon(h)=0,
\end{eqnarray*}
and the   antipode is defined by $S:H\rightarrow H,S(1_{H})=1_{H},S(h)=-h. $
\smallskip

 Recall from (\cite{CWZ2013}),  $A=k\{1_{A},x,g,gx\}$ is a Sweedler 4-dimension monoidal Hopf algebra  constructed from Sweedler 4-dimension Hopf algebra by Yau twist, where the twist map is defined by
 \begin{eqnarray*}
 \beta(1_{A})=1_{A}, \beta(g)=g, \beta(x)=-x, \beta(gx)=-gx,
\end{eqnarray*}
 the Hom-algebra structure $m$ is defined by
\begin{eqnarray*}
 && m(1_{A}\o1_{A})=1_{A}, m(1_{A}\o g)=g, m(1_{A}\o x)=-x, m(1_{A}\o gx)=-gx,\\
 && m(g\o 1_{A})=g, m(g\o g)=1, m(g\o x)=-gx, m(g\o gx)=-x,\\
 && m(x\o 1_{A})=-x, m(x\o g)=gx, m(x\o x)=0, m(x\o gx)=0,\\
  && m(gx\o 1_{A})=-gx, m(gx\o g)=x, m(gx\o x)=0, m(gx\o gx)=0,\\
\end{eqnarray*}
  the Hom-coalgebra structures $\epsilon$ and $\Delta$ are defined by
 \begin{eqnarray*}
 &&\epsilon(1_{A})=1,\epsilon(g)=\epsilon(x)=\epsilon(gx)=0,
 \Delta(1_{A})=1_{A}\o 1_{A}, \Delta(g)=g\o g,\\
 &&\Delta(x)=(-x)\o 1_{A}+g\o(-x),
\Delta(gx)=(-gx)\o g+1\o(-gx)
\end{eqnarray*}
and the   antipode is defined by $S:A\rightarrow A,S(1_{A})=1_{A},S(g)=g,S(x)=-gx,S(gx)=x. $

Now we define a left $(H,\a)$-Hom-module structure on $A$:
\begin{eqnarray*}
&&h\c 1_{A}=h\c g=h\c x=h\c gx=0,\\
 &&1_{H}\c 1_{A}=1_{A},1_{H}\c g=g,1_{H}\c x=-x,1_{H}\c gx=-gx.
\end{eqnarray*}
One may check directly that $A$ is a $(H,\a)$-Hom-module algebra.
Similarly,  we can define a left $(H,\a)$-Hom-comodule structure on $A$:
\begin{eqnarray*}
 \rho(1_{A})=1_{H}\o 1_{A},
 \rho(g)=1_{H}\o g,
 \rho(x)=1_{H}\o (-x),
 \rho(gx)=1_{H}\o (-gx).
\end{eqnarray*}
Then  $A$ is a $(H,\a)$-Hom-comodule algebra and $A$ is an object in $^{H}_{H}\mathcal{HYD}$.

Define the braiding  $C$  on $A$ by the usual flip map.
Clearly,  $C$ is symmetric on $A$.
By Proposition 2.2, there is a braided Hom-Lie algebra  $A^-$ with the bracket product [,] satisfying the following
non-vanishing relation
$$[x,g]=-[g,x]=2gx, [gx,g]=-[g,gx]=2x.$$

\noindent{\bf Lemma 2.6.}
Let $(A,\beta)$ be a monoidal Hom-algebra in $^{H}_{H}\mathcal{HYD}$
with monoidal Hom-subalgebras $X$ and $Y$
which are $H$-commutative such that $A=X+Y.$
Then the following equality holds:
\begin{eqnarray}
 &&\alpha^{-1}(u_{(-1)})\o\alpha^{-1}(y_{(-1)})\o(u_{0}y_{0})^{X}_{(-1)}\o(u_{0}y_{0})^{X}_{0}+\nonumber\\
&&\alpha^{-1}(u_{(-1)})\o\alpha^{-1}(y_{(-1)})\o(u_{0}y_{0})^{Y}_{(-1)}\o(u_{0}y_{0})^{Y}_{0}\nonumber\\
&=&u_{(-1)1}\o y_{(-1)1}\o u_{(-1)2}y_{(-1)2}\o\beta^{-1}((u_{0}y_{0})^{X})+\nonumber\\
&&u_{(-1)1}\o y_{(-1)1}\o u_{(-1)2}y_{(-1)2}\o\beta^{-1}((u_{0}y_{0})^{Y}),
\end{eqnarray}
for all $u,w\in X$ and $y,z\in Y$,
where $u_{0}y_{0}=(u_{0}y_{0})^{X}+(u_{0}y_{0})^{Y}\in X+Y.$
\medskip

{\bf Proof.} Since $\D(m_{(-1)})\o \beta^{-1}(m_{0})=\alpha^{-1}(m_{(-1)})\o (m_{0(-1)}\o m_{00})$,
by applying it to $u$ and $y$ respectively, we can get Eq. (2.4). $\hfill \Box$
\medskip

\noindent{\bf Lemma 2.7.}
Let $(A,\beta)$ be a monoidal Hom-algebra in $^{H}_{H}\mathcal{HYD}$
with monoidal Hom-subalgebras $X$ and $Y$
which are $H$-commutative such that $A=X+Y.$
Assume that the braiding $C$ is symmetric on $A$, then the following  equality holds:
\begin{eqnarray}
&&\epsilon(y_{(-1)})(\alpha(u_{(-1)})\c \beta^{-1}(w))\beta((u_{0}y_{0})^{X})-\nonumber\\
&&\epsilon(y_{(-1)})(\alpha(u_{(-1)})\c \beta^{-1}(z))\beta((u_{0}y_{0})^{Y})\nonumber\\
&=&\epsilon(u_{(-1)})\beta((u_{0}y_{0})^{X})(S^{-1}(\alpha(y_{(-1)}))\c\beta^{-1}(w))-\nonumber\\
&&\epsilon(u_{(-1)})\beta((u_{0}y_{0})^{Y})(S^{-1}(\alpha(y_{(-1)}))\c\beta^{-1}(z)),
\end{eqnarray}
for all $u,w\in X$ and $y,z\in Y$,
where $u_{0}y_{0}=(u_{0}y_{0})^{X}+(u_{0}y_{0})^{Y}\in X+Y.$
\medskip

{\bf Proof.}  For Eq. (2.5), we show it by the following computations:
\begin{eqnarray*}
 &&\epsilon(y_{(-1)})(\alpha(u_{(-1)})\c \beta^{-1}(w))\beta((u_{0}y_{0})^{X})-
\epsilon(y_{(-1)})(\alpha(u_{(-1)})\c \beta^{-1}(z))\beta((u_{0}y_{0})^{Y})\\
&=&\epsilon(y_{(-1)})((\alpha(u_{(-1)})\c \beta^{-1}(w))_{(-1)}\cdot(u_{0}y_{0})^{X})\beta((\alpha(u_{(-1)})\c \beta^{-1}(w))_{0})-\\
   &&\epsilon(y_{(-1)})((\alpha(u_{(-1)})\c \beta^{-1}(z))_{(-1)}\cdot(u_{0}y_{0})^{Y})\beta((\alpha(u_{(-1)})\c \beta^{-1}(z))_{0})\\
&=&\epsilon(y_{(-1)})\beta(\beta((u_{0}y_{0})^{X})_{0})(S^{-1}(\beta((u_{0}y_{0})^{X})_{(-1)})\cdot\beta^{-1}(\alpha(u_{(-1)})\c \beta^{-1}(w)))-\\
   &&\epsilon(y_{(-1)})\beta(\beta((u_{0}y_{0})^{Y})_{0})(S^{-1}(\beta((u_{0}y_{0})^{Y})_{(-1)})\cdot\beta^{-1}(\alpha(u_{(-1)})\c \beta^{-1}(z)))\\
&=&\epsilon(y_{(-1)})\beta^{2}((u_{0}y_{0})^{X}_{0})(S^{-1}(\alpha((u_{0}y_{0})^{X}_{(-1)}))\cdot\beta^{-1}(\alpha(u_{(-1)})\c \beta^{-1}(w)))-\\
   &&\epsilon(y_{(-1)})\beta^{2}((u_{0}y_{0})^{Y}_{0})(S^{-1}(\alpha((u_{0}y_{0})^{Y}_{(-1)}))\cdot\beta^{-1}(\alpha(u_{(-1)})\c \beta^{-1}(z))\\
&\stackrel{(2.4)}{=}&\epsilon(\alpha(y_{(-1)1}))\beta((u_{0}y_{0})^{X})(S^{-1}(\alpha(u_{(-1)2}y_{(-1)2}))\cdot\beta^{-1}(\alpha^{2}(u_{(-1)1})\c \beta^{-1}(w)))-\\
   &&\epsilon(\alpha(y_{(-1)1}))\beta((u_{0}y_{0})^{Y})(S^{-1}(\alpha(u_{(-1)2}y_{(-1)2}))\cdot\beta^{-1}(\alpha^{2}(u_{(-1)1})\c \beta^{-1}(z)))-\\
&=&\epsilon(u_{(-1)})(\beta((u_{0}y_{0})^{X})(S^{-1}(\alpha(y_{(-1)}))\cdot\beta^{-1}(w))-
     \beta((u_{0}y_{0})^{Y})(S^{-1}(\alpha(y_{(-1)}))\cdot\beta^{-1}(z))).
\end{eqnarray*}
The last equality holds since
\begin{eqnarray*}
&&\epsilon(\alpha(y_{(-1)1}))S^{-1}(\alpha(u_{(-1)2}y_{(-1)2}))\cdot\beta^{-1}(\alpha^{2}(u_{(-1)1})\c \beta^{-1}(w))\\
&=&\epsilon(\alpha(y_{(-1)1}))S^{-1}(\alpha(u_{(-1)2}y_{(-1)2}))\cdot(\alpha(u_{(-1)1})\c \beta^{-2}(w))\\
&=&\epsilon(y_{(-1)1})((S^{-1}(y_{(-1)2})S^{-1}(u_{(-1)2}))\alpha(u_{(-1)1}))\c \beta^{-1}(w)\\
&=&\epsilon(y_{(-1)1})(\alpha(S^{-1}(y_{(-1)2}))(S^{-1}(u_{(-1)2})u_{(-1)1}))\c \beta^{-1}(w)\\
&=&(S^{-1}(y_{(-1)})(\epsilon(u_{(-1)})1_{H}))\c \beta^{-1}(w)\\
&=&\epsilon(u_{(-1)})S^{-1}(\alpha(y_{(-1)}))\c \beta^{-1}(w).
\end{eqnarray*}
And this completes the proof.
$\hfill \Box$
\medskip

\noindent{\bf Theorem 2.8.}
Let $(A,\beta)$ be a monoidal Hom-algebra in $^{H}_{H}\mathcal{HYD}$
with monoidal Hom-subalgebras $X$ and $Y$
which are $H$-commutative such that $A=X+Y.$
Assume that the braiding $C$ is symmetric on $A$, then $[A,A][A,A]=0.$
\medskip

{\bf Proof.} It is sufficient to prove $[u,x][v,y]=0$ holds for all $u,v\in X$
and $x,y\in Y$.
For any $a,b,c,d\in A$, we first note that
$(ab)(cd)=(a\beta^{-1}(bc))\beta(d)$ which can be verified easily from the Hom-associativity of $A$.
By the definition of the bracket product, we have
\begin{eqnarray*}
[u,x][v,y]
&=&(ux-(u_{(-1)}\c \beta^{-1}(x))\beta(u_{0}))(vy-(v_{(-1)}\c \beta^{-1}(y))\beta(v_{0}))\\
&=&(ux)(vy)+((u_{(-1)}\c \beta^{-1}(x))\beta(u_{0}))((v_{(-1)}\c \beta^{-1}(y))\beta(v_{0}))-\\
&&(ux)((v_{(-1)}\c \beta^{-1}(y))\beta(v_{0}))-((u_{(-1)}\c \beta^{-1}(x))\beta(u_{0}))(vy).
\end{eqnarray*}

Next we will compute the four expressions above respectively.
For this purpose,
let $xv=w+z,$ where $w\in X, z\in Y$.

(1)
$
(ux)(vy)=((u_{(-1)}\cdot\beta^{-2}(w))\beta(u_{0}))\beta(y)
+(u\beta^{-1}(z_{(-1)}\cdot y))\beta(z_{0}).
$
In  fact,
\begin{eqnarray*}
(ux)(vy)&=&(u\beta^{-1}(xv))\beta(y)
=(u\beta^{-1}(w))\beta(y)+\beta(u)(\beta^{-1}(z)y)\\
&=&((u_{(-1)}\cdot\beta^{-2}(w))\beta(u_{0}))\beta(y)+\beta(u)((\alpha^{-1}(z_{(-1)})\cdot\beta^{-1}(y))\beta(\beta^{-1}(z_{0})))\\
&=&((u_{(-1)}\cdot\beta^{-2}(w))\beta(u_{0}))\beta(y)+\beta(u)((\alpha^{-1}(z_{(-1)})\cdot\beta^{-1}(y))z_{0})\\
&=&((u_{(-1)}\cdot\beta^{-2}(w))\beta(u_{0}))\beta(y)
+(u\beta^{-1}(z_{(-1)}\cdot y))\beta(z_{0}).
\end{eqnarray*}

(2) $((u_{(-1)}\c \beta^{-1}(x))\beta(u_{0}))(vy)
=((u_{(-1)}\cdot\beta^{-2}(w))\beta(u_{0}))\beta(y)
+\epsilon(y_{(-1)})(\alpha(u_{(-1)})\cdot\beta^{-1}(z))\beta((u_{0}y_{0})^{X})
+\epsilon(y_{(-1)})(\alpha(u_{(-1)})\cdot\beta^{-1}(z))\beta((u_{0}y_{0})^{Y})$.
In  fact,
\begin{eqnarray*}
&&((u_{(-1)}\c \beta^{-1}(x))\beta(u_{0}))(vy)\\
&=&((u_{(-1)}\c \beta^{-1}(x))\beta^{-1}(\beta(u_{0})v))\beta(y)\\
&=&((u_{(-1)}\c \beta^{-1}(x))(u_{0}\beta^{-1}(v)))\beta(y)\\
&=&((u_{(-1)}\c \beta^{-1}(x))((u_{0(-1)}\cdot\beta^{-2}(v))\beta(u_{00})))\beta(y)\\
&=&((\alpha(u_{(-1)1})\c \beta^{-1}(x))((u_{(-1)2}\cdot\beta^{-2}(v))u_{0}))\beta(y)\\
&=&(((u_{(-1)1}\c \beta^{-2}(x))(u_{(-1)2}\cdot\beta^{-2}(v)))\beta(u_{0}))\beta(y)\\
&=&((u_{(-1)}\c \beta^{-2}(xv))\beta(u_{0}))\beta(y)\\
&=&((u_{(-1)}\c \beta^{-2}(w))\beta(u_{0}))\beta(y)+((u_{(-1)}\c \beta^{-2}(z))\beta(u_{0}))\beta(y)\\
&=&((u_{(-1)}\c \beta^{-2}(w))\beta(u_{0}))\beta(y)+(\alpha(u_{(-1)})\c \beta^{-1}(z))(\beta(u_{0})\beta(y_{0}))\epsilon(y_{(-1)})\\
&=&((u_{(-1)}\c \beta^{-2}(w))\beta(u_{0}))\beta(y)+(\alpha(u_{(-1)})\c \beta^{-1}(z))\beta(u_{0}y_{0})\epsilon(y_{(-1)})\\
&=&((u_{(-1)}\c \beta^{-2}(w))\beta(u_{0}))\beta(y)+\epsilon(y_{(-1)})(\alpha(u_{(-1)})\c \beta^{-1}(z))\beta((u_{0}y_{0})^{X})\\
    &&+\epsilon(y_{(-1)})(\alpha(u_{(-1)})\c \beta^{-1}(z))\beta((u_{0}y_{0})^{Y}).
\end{eqnarray*}

(3) $(ux)((v_{(-1)}\c \beta^{-1}(y))\beta(v_{0}))=(u \beta^{-1}(z_{(-1)}\c y))\beta(z_{0})+\epsilon(u_{(-1)})\beta((u_{0}y_{0})^{X})(S^{-1}(\alpha(y_{(-1)}))\cdot\beta^{-1}(w))+
   \epsilon(u_{(-1)})\beta((u_{0}y_{0})^{Y})(S^{-1}(\alpha(y_{(-1)}))\cdot\beta^{-1}(w))$. In fact,
\begin{eqnarray*}
&&(ux)((v_{(-1)}\c \beta^{-1}(y))\beta(v_{0}))\\
&=&(u \beta^{-1}(x(v_{(-1)}\c \beta^{-1}(y))))\beta^{2}(v_{0})\\
&=&(u \beta^{-1}((x_{(-1)}\cdot\beta^{-1}(v_{(-1)}\c \beta^{-1}(y)))\beta(x_{0})))\beta^{2}(v_{0})\\
&=&(u \beta^{-1}((x_{(-1)}\cdot(\alpha^{-1}(v_{(-1)})\c \beta^{-2}(y)))\beta(x_{0})))\beta^{2}(v_{0})\\
&=&(u \beta^{-1}((\alpha^{-1}(x_{(-1)}v_{(-1)})\c \beta^{-1}(y))\beta(x_{0}))\beta^{2}(v_{0})\\
&=&\beta(u)(((\alpha^{-2}(x_{(-1)}v_{(-1)})\c \beta^{-2}(y))x_{0})\beta(v_{0}))\\
&=&\beta(u)((\alpha^{-1}(x_{(-1)}v_{(-1)})\c \beta^{-1}(y))(x_{0}v_{0}))\\
&=&(u \beta^{-1}((x_{(-1)}v_{(-1)})\c y))\beta(x_{0}v_{0})\\
&=&(u \beta^{-1}((xv)_{(-1)}\c y))\beta((xv)_{0})\\
&=&(u \beta^{-1}(w_{(-1)}\c y))\beta(w_{0})+(u \beta^{-1}(z_{(-1)}\c y))\beta(z_{0})\\
&=&(u \beta(y_{0}))(S^{-1}(\alpha(y_{(-1)}))\cdot\beta^{-1}(w))+(u \beta^{-1}(z_{(-1)}\c y))\beta(z_{0})
 \end{eqnarray*}
   \begin{eqnarray*}
&=&\epsilon(u_{(-1)})\beta(u_{0}y_{0})(S^{-1}(\alpha(y_{(-1)}))\cdot\beta^{-1}(w))+(u \beta^{-1}(z_{(-1)}\c y))\beta(z_{0})\\
&=&\epsilon(u_{(-1)})\beta((u_{0}y_{0})^{X})(S^{-1}(\alpha(y_{(-1)}))\cdot\beta^{-1}(w))+(u \beta^{-1}(z_{(-1)}\c y))\beta(z_{0})+\\
   &&\epsilon(u_{(-1)})\beta((u_{0}y_{0})^{Y})(S^{-1}(\alpha(y_{(-1)}))\cdot\beta^{-1}(w)).
\end{eqnarray*}

(4) $((u_{(-1)}\c \beta^{-1}(x))\beta(u_{0}))((v_{(-1)}\c \beta^{-1}(y))\beta(v_{0}))=\epsilon(y_{(-1)})(\alpha(u_{(-1)})\c \beta^{-1}(w))\beta((u_{0}y_{0})^{X})\\+
\epsilon(y_{(-1)})(\alpha(u_{(-1)})\c \beta^{-1}(w))\beta((u_{0}y_{0})^{X})+
  \epsilon(u_{(-1)})\beta((u_{0}y_{0})^{Y})(S^{-1}(\alpha(y_{(-1)}))\cdot \beta^{-1}(w))+\\
  \epsilon(u_{(-1)})\beta((u_{0}y_{0})^{Y})(S^{-1}(\alpha(y_{(-1)}))\cdot \beta^{-1}(z))$.

 Here we first give two useful equalities:
 \begin{eqnarray}
(u_{(-1)2}y_{(-1)2})\cdot(S^{-1}(y_{(-1)1})\cdot\beta^{-2}(v))
 &=&\epsilon(y_{(-1)})\alpha(u_{(-1)2})\cdot\beta^{-1}(v)),\\
 (S^{-1}(y_{(-1)2})S^{-1}(u_{(-1)2}))\cdot(u_{(-1)1}\cdot\beta^{-2}(v))
 &=&\epsilon(u_{(-1)})S^{-1}(\alpha(y_{(-1)2}))\cdot\beta^{-1}(v).~~~~~~~
 \end{eqnarray}
 In fact,
  \begin{eqnarray*}
&&(u_{(-1)2}y_{(-1)2})\cdot(S^{-1}(y_{(-1)1})\cdot\beta^{-2}(v))\\
&=&((\alpha^{-1}(u_{(-1)2})\alpha^{-1}(y_{(-1)2}))S^{-1}(y_{(-1)1}))\cdot\beta^{-1}(v)\\
&=&(u_{(-1)2}(\alpha^{-1}(y_{(-1)2})\alpha^{-1}(S^{-1}(y_{(-1)1}))))\cdot\beta^{-1}(v)\\
&=&(u_{(-1)2}\epsilon(y_{(-1)})1_{H})\cdot\beta^{-1}(v)
=\epsilon(y_{(-1)})\alpha(u_{(-1)2})\cdot\beta^{-1}(v).
 \end{eqnarray*}
So Eq. (2.6) holds and similarly for Eq. (2.7).
 Therefore,
\begin{eqnarray*}
&&((u_{(-1)}\c \beta^{-1}(x))\beta(u_{0}))((v_{(-1)}\c \beta^{-1}(y))\beta(v_{0}))\\
&=&((u_{(-1)}\c \beta^{-1}(x))\beta(u_{0}))(\beta(y_{0})(S^{-1}(y_{(-1)})\cdot \beta^{-1}(v)))\\
&=&((u_{(-1)}\c \beta^{-1}(x))(u_{0}y_{0}))\beta(S^{-1}(y_{(-1)})\cdot \beta^{-1}(v))\\
&=&((u_{(-1)}\c \beta^{-1}(x))(u_{0}y_{0}))(S^{-1}(\alpha(y_{(-1)}))\cdot v)\\
&=&\beta(u_{(-1)}\c \beta^{-1}(x))((u_{0}y_{0})(S^{-1}(y_{(-1)})\cdot \beta^{-1}(v)))\\
&=&\beta(u_{(-1)}\c \beta^{-1}(x))((u_{0}y_{0})^{X}(S^{-1}(y_{(-1)})\cdot \beta^{-1}(v)))+\\
   &&\beta(u_{(-1)}\c \beta^{-1}(x))((u_{0}y_{0})^{Y}(S^{-1}(y_{(-1)})\cdot \beta^{-1}(v)))\\
&=&\beta(u_{(-1)}\c \beta^{-1}(x))(((u_{0}y_{0})^{X}_{(-1)}\cdot\beta^{-1}(S^{-1}(y_{(-1)})\cdot \beta^{-1}(v)))\beta((u_{0}y_{0})^{X}_{0}))+\\
   &&((u_{(-1)}\c \beta^{-1}(x))(u_{0}y_{0})^{Y})\beta(S^{-1}(y_{(-1)})\cdot \beta^{-1}(v))\\
&=&((u_{(-1)}\c \beta^{-1}(x))((u_{0}y_{0})^{X}_{(-1)}\cdot\beta^{-1}(S^{-1}(y_{(-1)})\cdot \beta^{-1}(v))))\beta^{2}((u_{0}y_{0})^{X}_{0})+\\
   &&(((u_{(-1)}\c \beta^{-1}(x))_{(-1)}\cdot\beta^{-1}((u_{0}y_{0})^{Y}))\beta((u_{(-1)}\c \beta^{-1}(x))_{0}))\beta(S^{-1}(y_{(-1)})\cdot \beta^{-1}(v))\\
&=&((u_{(-1)}\c \beta^{-1}(x))((u_{0}y_{0})^{X}_{(-1)}\cdot\beta^{-1}(S^{-1}(y_{(-1)})\cdot \beta^{-1}(v))))\beta^{2}((u_{0}y_{0})^{X}_{0})+\\
   &&(\beta((u_{0}y_{0})^{Y}_{0})(S^{-1}((u_{0}y_{0})^{Y}_{(-1)})\cdot\beta^{-1}(u_{(-1)}\c \beta^{-1}(x))))\beta(S^{-1}(y_{(-1)})\cdot \beta^{-1}(v))\\
&=&((\alpha(u_{(-1)1})\c \beta^{-1}(x))((u_{(-1)2}y_{(-1)2})\cdot\beta^{-1}(S^{-1}(\alpha(y_{(-1)1}))\cdot \beta^{-1}(v))))\beta((u_{0}y_{0})^{X})+\\
   &&((u_{0}y_{0})^{Y}(S^{-1}(u_{(-1)2}y_{(-1)2})\cdot\beta^{-1}(\alpha(u_{(-1)1})\c \beta^{-1}(x))))\beta(S^{-1}(\alpha(y_{(-1)1}))\cdot \beta^{-1}(v))\\
&=&((\alpha(u_{(-1)1})\c \beta^{-1}(x))((u_{(-1)2}y_{(-1)2})\cdot(S^{-1}(y_{(-1)1})\cdot \beta^{-2}(v))))\beta((u_{0}y_{0})^{X})+\\
   &&((u_{0}y_{0})^{Y}((S^{-1}(u_{(-1)2})S^{-1}(y_{(-1)2}))\cdot(u_{(-1)1}\c \beta^{-2}(x))))\beta(S^{-1}(\alpha(y_{(-1)1}))\cdot \beta^{-1}(v))
    \end{eqnarray*}
   \begin{eqnarray*}
&\stackrel{(2.6),(2.7)}{=}&\epsilon(y_{(-1)})((\alpha(u_{(-1)1})\c \beta^{-1}(x))(\alpha(u_{(-1)2})\beta^{-1}(v)))\beta((u_{0}y_{0})^{X})+\\
  &&\epsilon(u_{(-1)})((u_{0}y_{0})^{Y}(S^{-1}(\alpha(y_{(-1)2}))\cdot \beta^{-1}(x)))\beta(S^{-1}(\alpha(y_{(-1)1}))\cdot \beta^{-1}(v))\\
&=&\epsilon(y_{(-1)})(\alpha(u_{(-1)})\c \beta^{-1}(xv))\beta((u_{0}y_{0})^{X})+
   \epsilon(u_{(-1)})\beta((u_{0}y_{0})^{Y})(S^{-1}(\alpha(y_{(-1)}))\cdot \beta^{-1}(xv))\\
&=&\epsilon(y_{(-1)})(\alpha(u_{(-1)})\c \beta^{-1}(w))\beta((u_{0}y_{0})^{X})+
\epsilon(y_{(-1)})(\alpha(u_{(-1)})\c \beta^{-1}(w))\beta((u_{0}y_{0})^{X})+\\
  &&\epsilon(u_{(-1)})\beta((u_{0}y_{0})^{Y})(S^{-1}(\alpha(y_{(-1)}))\cdot \beta^{-1}(w))+
  \epsilon(u_{(-1)})\beta((u_{0}y_{0})^{Y})(S^{-1}(\alpha(y_{(-1)}))\cdot \beta^{-1}(z)).
\end{eqnarray*}

Hence we have
\begin{eqnarray*}
[u,x][v,y]
&=&-\epsilon(y_{(-1)})(\alpha(u_{(-1)})\c \beta^{-1}(z))\beta((u_{0}y_{0})^{Y})\\
&&-\epsilon(u_{(-1)})\beta((u_{0}y_{0})^{X})(S^{-1}(\alpha(y_{(-1)}))\c\beta^{-1}(w))\\
&&+\epsilon(y_{(-1)})(\alpha(u_{(-1)})\c\beta^{-1}(w))\beta((u_{0}y_{0})^{X})\\
&&+\epsilon(u_{(-1)})\beta((u_{0}y_{0})^{Y})(S^{-1}(\alpha(y_{(-1)}))\c \beta^{-1}(z))\\
&\stackrel{(3.2)}{=}&0,
\end{eqnarray*}
as desired. And this completes the proof.$\hfill \Box$
\medskip

\noindent{\bf Corollary 2.9.}
Under the hypotheses of the theorem above, $[A,A]$ is nilpotent.
If $A$ is also $H$-semiprime, then $A$ is $H$-commutative.

{\bf Proof}  Straightforward from Theorem 2.8.$\hfill \Box$
\medskip

\section{Central invariants of braided Hom-Lie algebras}
\def\theequation{\arabic{section}.\arabic{equation}}
\setcounter{equation} {0}

In this section, we always assume that $(H,\alpha)$ is a monoidal Hom-Hopf algebra.
We consider some $H$-analogous of classical concepts of ring theory
and of Lie theory as follows.

Let $(A,\beta)$ be be a monoidal Hom-algebra in $^{H}_{H}\mathcal{HYD}$.
An {\em $H$-Hom-ideal} $U$ of $A$ is not only $H$-stable but also $H$-costable
such that $\beta(U)\subseteq U$ and $(AU)A=A(UA)\subseteq U.$

Let $(L,\beta)$ be a braided Hom-Lie algebra.
An {\em $H$-Hom-Lie ideal} $U$ of $L$ is not only $H$-stable but also $H$-costable
such that $\beta(U)\subseteq U$ and $[U,L]\subseteq U.$

Define the {\em center} of $L$ to be $Z_{H}(L)=\{l\in L|[l,L]_{H}=0\}.$
It is easy to see that $Z_{H}(L)$ is not only $H$-stable but also $H$-costable.

$L$ is called {\em $H$-prime} if the product of any two non-zero $H$-Hom-ideals of $L$ is non-zero.
It is called {\em $H$-semiprime} if it has no  non-zero nilpotent $H$-Hom-ideals,
and is called {\em $H$-simple} if it has no nontrivial $H$-Hom-ideals.

\noindent{\bf Definition 3.1.}
If $(A,\beta)$ is a monoidal Hom-algebra in $^{H}_{H}\mathcal{HYD}$, the monoidal Hom-subalgebra
of $H$-invariant is the set:
$$
A_0=\{a\in A|h\c a=\v(h)a,~\mbox{for all}~h\in H\}.
$$

Recall from Proposition 2.2, a monoidal Hom-algebra $(L,\beta)$ in $^{H}_{H}\mathcal{HYD}$
gives rise to a braided Hom-Lie algebra $(L,[,],\beta)$ in $^{H}_{H}\mathcal{HYD}$.

In what follows, we always assume that the bracket product in braided Hom-Lie algebra $(L,[,],\beta)$ is defined as Proposition 2.2, that is .
\begin{eqnarray}
[,]: A\otimes A\rightarrow A~by~[a,b]=ab-(a_{(-1)}\cdot \beta^{-1}(b))\beta(a_{0}),~a,b\in A.
\end{eqnarray}

\noindent{\bf Lemma 3.2.}
Let $(L,\beta)$ be a monoidal Hom-algebra in $^{H}_{H}\mathcal{HYD}$
 and $(L,[,],\beta)$ the derived braided Hom-Lie algebra. Then

(1) $[\beta(a),bc]=[a,b]\beta(c)+(\alpha(a_{(-1)})\c b)[\beta(a_{0}),c],$

(2) $[ab,\beta(c)]=\beta(a)[b,c]+[a,b_{(-1)}\c \beta^{-1}(c)]\beta^{2}(b_{0}),$
for all $a,b,c\in L.$
\medskip

{\bf Proof.}
(1) For all $a,b,c\in L$, it is clear that
 $[a,b]\beta(c)=(ab)\beta(c)-((a_{(-1)}\cdot \beta^{-1}(b))\beta(a_{0}))\beta(c)$.
Similarly,
\begin{eqnarray*}
&&(\alpha(a_{(-1)})\c b)[\beta(a_{0}),c]\\
&=&(\alpha(a_{(-1)})\c b)(\beta(a_{0})c)-(\alpha(a_{(-1)})\c b)((\alpha(a_{0(-1)})\cdot \beta^{-1}(c))\beta^{2}(a_{00}))\\
&=&\beta(a_{(-1)}\c\beta^{-1}(b))(\beta(a_{0})c)-\beta(a_{(-1)}\c \beta^{-1}(b))((\alpha(a_{0(-1)})\cdot \beta^{-1}(c))\beta^{2}(a_{00}))\\
&=&((a_{(-1)}\c\beta^{-1}(b))\beta(a_{0}))\beta(c)-((a_{(-1)}\c \beta^{-1}(b))(\alpha(a_{0(-1)})\cdot \beta^{-1}(c)))\beta^{3}(a_{00}))\\
&=&((a_{(-1)}\c\beta^{-1}(b))\beta(a_{0}))\beta(c)-((\alpha(a_{(-1)1})\c \beta^{-1}(b))(\alpha(a_{(-1)2})\cdot \beta^{-1}(c)))\beta^{2}(a_{0}))\\
&=&((a_{(-1)}\c\beta^{-1}(b))\beta(a_{0}))\beta(c)-(\alpha(a_{(-1)})\c \beta^{-1}(bc))\beta^{2}(a_{0})).
\end{eqnarray*}
Therefore,
\begin{eqnarray*}
&&[a,b]\beta(c)+(\alpha(a_{(-1)})\c b)[\beta(a_{0}),c]\\
&=&(ab)\beta(c)-(\alpha(a_{(-1)})\c \beta^{-1}(bc))\beta^{2}(a_{0}))\\
&=&\beta(a)(bc)-((\alpha(a_{(-1)})\cdot \beta^{-1}(bc))\beta^{2}(a_{0})\\
&=&\beta(a)(bc)-((\beta(a))_{(-1)}\cdot \beta^{-1}(bc))\beta((\beta(a))_{0})\\
&=&[\beta(a),bc].
\end{eqnarray*}

(2) For all $a,b,c\in L$, on the one hand, we have
\begin{eqnarray*}
\beta(a)[b,c]
&=&\beta(a)(bc)-\beta(a)((b_{(-1)}\cdot \beta^{-1}(c))\beta(b_{0}))\\
&=&(ab)\beta(c)-(a(b_{(-1)}\cdot \beta^{-1}(c)))\beta^{2}(b_{0}).
\end{eqnarray*}
On the other hand, we get
\begin{eqnarray*}
&&[a,b_{(-1)}\c \beta^{-1}(c)]\beta^{2}(b_{0})\\
&=&(a(b_{(-1)}\c \beta^{-1}(c)))\beta^{2}(b_{0})-((a_{(-1)}\cdot\beta^{-1}(b_{(-1)}\c \beta^{-1}(c)))\beta(a_{0}))\beta^{2}(b_{0})\\
&=&(a(b_{(-1)}\c \beta^{-1}(c)))\beta^{2}(b_{0})-((a_{(-1)}\cdot(\alpha^{-1}(b_{(-1)})\c \beta^{-2}(c)))\beta(a_{0}))\beta^{2}(b_{0})\\
&=&(a(b_{(-1)}\c \beta^{-1}(c)))\beta^{2}(b_{0})-(((\alpha^{-1}(a_{(-1)})\alpha^{-1}(b_{(-1)}))\c \beta^{-1}(c))\beta(a_{0}))\beta^{2}(b_{0})\\
&=&(a(b_{(-1)}\c \beta^{-1}(c)))\beta^{2}(b_{0})-(a_{(-1)}b_{(-1)}\c c)\beta(a_{0}b_{0}).
\end{eqnarray*}
It follows that
\begin{eqnarray*}
&&\beta(a)[b,c]+[a,b_{(-1)}\c \beta^{-1}(c)]\beta^{2}(b_{0})\\
&=&\beta(a)(bc)-(a_{(-1)}b_{(-1)}\c c)\beta(a_{0}b_{0})\\
&=&(ab)\beta(c)-(a_{(-1)}b_{(-1)}\c c)\beta(a_{0}b_{0})\\
&=&[ab,\beta(c)].
\end{eqnarray*}
The proof is completed.$\hfill \Box$
\medskip

Define $ad_x(l)=[x,l]$ for all $x, l\in L$, By Lemma 3.2(1) we have
$$
ad_{\beta(x)}(lm)=ad_x(l)\alpha(m)+(\alpha^{-1}(x_{(-1)})\cdot \beta(l))ad_{x_0}(m),~x, l, m\in L.
$$

\noindent{\bf Lemma 3.3.}
Let $(L,\beta)$ be a monoidal Hom-algebra in $^{H}_{H}\mathcal{HYD}$ and $x$ an $\beta$-invariant element in $L_0$.
 Then for any $y, z\in L$, the following equalities hold:

(1) $C_{L,L}(x\o y)=y\o x$, $C_{L,L}(y\o x)=x\o y$;

(2) $ad_x(y)=xy-yx$;

(3) $ad_x(yz)=ad_x(y)\beta(z)+\beta(y)ad_x(z)$;

(4) $ad^2_x(yz)=ad^2_x(y)\beta^2(z)+2\beta(ad_x(y)ad_x(z))+\beta^2(y)ad^2_x(z)$.
\medskip

{\bf Proof.}
(1) Since $x\in L_0$, we have
\begin{eqnarray*}
C_{L,L}(y\o x)
&=&y_{(-1)}\cdot \beta^{-1}(x)\otimes\beta(y_{0})
=y_{(-1)}\cdot x\otimes\beta(y_{0})\\
&=&\epsilon(y_{(-1)})x\otimes\beta(y_{0})
=x\o y,\\
C_{L,L}(x\o y)
&=&x_{(-1)}\cdot \beta^{-1}(y)\otimes\beta(x_{0})
=\beta(y_{0})\otimes S^{-1}(y_{(-1)})\cdot \beta^{-1}(x)\\
&=&\beta(y_{0})\otimes S^{-1}(y_{(-1)})\cdot x
=\beta(y_{0})\otimes \epsilon(S^{-1}(y_{(-1)}))x
=y\o x.
\end{eqnarray*}

(2) Straightforward from (1).

(3) Straightforward from Lemma 3.2 (1).

(4) By (2) and (3), we have
\begin{eqnarray*}
ad^2_x(yz)
&=&ad_x(ad_x(y)\beta(z)+\beta(y)ad_x(z))\\
&=&ad_x(ad_x(y)\beta(z))+ad_x(\beta(y)ad_x(z))\\
&=&ad^2_x(y)\beta^2(z)+\beta(ad_x(y))ad_x\beta(z)+\\
   &&ad_x\beta(y)\beta(ad_x(z))+\beta^2(y)ad^2_x(z)\\
&=&ad^2_x(y)\beta^2(z)+\beta(ad_x(y))ad_{\beta(x)}\beta(z)+\\
   &&ad_{\beta(x)}\beta(y)\beta(ad_x(z))+\beta^2(y)ad^2_x(z)\\
&=&ad^2_x(y)\beta^2(z)+2\beta(ad_x(y)ad_x(z))+\beta^2(y)ad^2_x(z).
\end{eqnarray*}
The proof is finished.
$\hfill \Box$
\medskip

\noindent{\bf Lemma 3.4.}
Let $(L,[,],\beta)$ be the derived braided Hom-Lie algebra.
 Assume that $L$ is $H$-simple, then $Z_H(L)_0$ is a field.
\medskip

{\bf Proof.}
Note that  $Z_H(L)_0=Z_H(L)\cap L_0=Z(L)\cap L_0=Z(L)_0$, where $Z(L)$ is the
usual center of $L$. Taking $0\neq x\in Z_H(L)_0$, we have that $Lx=I\neq 0$ is an $H$-Hom-ideal,
thus $I=L$ since $L$ is $H$-simple. That is to say that for some $y\in L$, we obtain $xy=yx=1$. Since
\begin{eqnarray*}
\beta^2(h\cdot y)
&=&\beta(h\cdot y)1
=\beta(h\cdot y)(xy)\\
&=&\beta(\alpha(h_1)\cdot y)(\epsilon(\alpha(h_2))xy)\\
&=&\beta(\alpha(h_1)\cdot y)((\alpha(h_2)\cdot x)y)\\
&=&((\alpha(h_1)\cdot y)(\alpha(h_2)\cdot x))\beta(y)\\
&=&(\alpha(h)\cdot (xy))\beta(y)
=(\alpha(h)\cdot 1)\beta(y)\\
&=&(\epsilon(\alpha(h))1)\beta(y)
=\epsilon(h)\beta^{2}(y)\\
&=&\beta^{2}(\epsilon(h)y)
\end{eqnarray*}
We can get $h\cdot y=\v(h)y$ since $\beta$ is bijective, that is, $y\in L_0$.

We need to show $y\in Z_H(L)$. For any $z\in L$, by Lemma 3.3(1), $[z,x]=zx-xz=0$.
Then we  have
\begin{eqnarray*}
&&\beta^2(yz-zy)
=\beta^2(yz)-\beta^2(zy)\\
&=&\beta(yz)\beta(1)-\beta(yx)\beta(zy)\\
&=&\beta^{2}(y)(\beta(z)1)-\beta^{2}(y)(\beta(x)(zy))\\
&=&\beta^{2}(y)(\beta(z)(xy))-\beta^{2}(y)(\beta(x)(zy))\\
&=&\beta^{2}(y)((zx)\beta(y))-\beta^{2}(y)((xz)\beta(y))\\
&=&\beta^{2}(y)((zx-xz)\beta(y))\\
&=&0.
\end{eqnarray*}
Since $\beta$ is bijective, it follows that $yz=zy$,
 i.e. $[y,z]=yz-zy=0$ by Lemma 3.3 (2). This shows that $y\in Z_H(L)$, as desired.
$\hfill \Box$
\medskip

\noindent{\bf Lemma 3.5.}
Let $(L, [,], \b)$  be the derived braided Hom-Lie algebra and $x$ an $\beta$-invariant element in $L_0$, $l, m\in L$. Then

(1) $ad^2_x(xl)=xad^2_x(l)$;

(2) If $ad^2_x(L)=0$ and char$(k)\neq 2$, then $ad_x(l)(Lad_x(m))=0$.
\medskip

{\bf Proof.}
(1) It is straightforward from Lemma 3.3 (4).

(2) For all $l, m\in L$, we have
\begin{eqnarray*}
0&=&ad^2_x(lm)=ad^2_x(l)\beta^2(m)+2\beta(ad_x(l)ad_x(m))+\beta^2(l)ad^2_x(m)\\
&=&2ad_x(\beta(l))ad_x(\beta(m)).
\end{eqnarray*}
So $ad_x(l)ad_x(m)=0$ since char$(k)\neq 2$.
For any $z\in L$, by Lemma 3.3 (3),
$z ad_x(m)=ad_x(\beta^{-1}(z)m)-ad_x(\beta^{-1}(z))\beta(m)$.
Therefore,
\begin{eqnarray*}
ad_x(l)(z ad_x(m))
&=&ad_x(l)ad_x(\beta^{-1}(z)m)-ad_x(l)(ad_x(\beta^{-1}(z))\beta(m))\\
&=&0-\beta(ad_x(\beta^{-1}(l)))(ad_x(\beta^{-1}(l))\beta(m))\\
&=&-(ad_x(\beta^{-1}(l))ad_x(\beta^{-1}(l)))m\\
&=&0.
\end{eqnarray*}
By the arbitrary of $z$, $ad_x(l)(Lad_x(m))=0$. And this finishes the proof.$\hfill \Box$
\medskip

\noindent{\bf Lemma 3.6.}
Let $(L, [,], \beta)$ be the derived braided Hom-Lie algebra and  $I$ an $H$-Hom-Lie ideal of $[L,L]$.
Assume that $L$ is $H$-simple and char$(k)\neq 2$.
If $x$ is an $\beta$-invariant element in $I_0$ satisfying (i) $ad_x(I)=0$, (ii) $ad^2_x([L,L])=0$.
 Then $x\in Z_H(L)$.
\medskip

{\bf Proof.}
For any $m\in L$, $l\in [L,L]$ and $y\in I$. By Lemma 3.2 (1),
\begin{eqnarray*}
0=ad^2_x([\beta(l), my])=ad^2_x([l, m]\beta(y))+ad^2_x((\alpha(l_{(-1)})\cdot m)[\beta(l_0), y]).
\end{eqnarray*}
 First, we have
\begin{eqnarray*}
&& ad^2_x([l, m]\beta(y))\\
&=& ad^2_x([l, m])\beta^3(y)+2\beta(ad_x([l, m])ad_x(\beta(y)))+\beta^2([l,m])ad^2_x(\beta(y))\\
&\stackrel{(i)}{=}& ad^2_x([l, m])\beta^3(y)
\stackrel{(ii)}{=}0.
\end{eqnarray*}
So $ad^2_x((\alpha(l_{(-1)})\cdot m)[\beta(l_0), y])$.
On the other hand, since $l\in [L,L]$ and $[,]$ is $H$-colinear, it follows that $\beta(l_0)\in [L,L]$,
$ad_x([l_0, y])\stackrel{(i)}{=}0$ and $ad^2_x([l_0, y])\stackrel{(ii)}{=}0$.
Therefore,
\begin{eqnarray*}
&& ad^2_x(\alpha(l_{(-1)})\cdot m)[\beta(l_0), y])\\
&=& ad^2_x(\alpha(l_{(-1)})\cdot m)\beta^2([\beta(l_0), y])
+ 2\beta(ad_x(\alpha(l_{(-1)})\cdot m)ad_x([\beta(l_0), y]))\\
&&+\beta^2(\alpha(l_{(-1)})\cdot m)ad^2_x([\beta(l_0), y])\\
&=&ad^2_x(\alpha(l_{(-1)})\cdot m)\beta^2([\beta(l_0), y]).
\end{eqnarray*}
Thus we obtain
$
ad^2_x(\alpha(l_{(-1)})\cdot m)\beta^2([\beta(l_0), y])=0.
$
We completes the proof by the following two cases:

Case (1): If $[I,[L,L]]=0$, then we have $ad^2_x(L)=0$. By Lemma 3.5 (2), $ad_x(l)(Lad_x(m))=0$.
Since $L$ is $H$-simple, we get $ad_x(l)=0$. So $x\in Z_H(L)$ since $l$ is an arbitrary element in $L$.

Case (2): If $[I,[L,L]]\neq 0$, let $U=[I,[L,L]]$.
It is easy to see that $U$ is a $H$-Hom-Lie ideal of $[L,L]$.
Since $ad^2_x(\alpha(l_{(-1)})\cdot m)\beta^2([\beta(l_0), y])=0$, we have $ad^2_x(L)U=0$.
Let $Q=\{y\in L|yU=0\}$, then $Q$ is an $H$-stable
$H$-costable left Hom-ideal of $L$, we claim $Q=0$. If not, then $L=QL$ since $L$ is $H$-simple. By
(2.1) we have
$$
QL\subseteq [Q,L]+LQ\subseteq [Q,L]+Q\subseteq L.
$$
Thus $L=Q+[Q,L]$. Let $y\in Q$, $l\in [L,L]$ and $u\in U$. Since $Q$ is an $H$-Hom-ideal, $\beta^{2}(y_{0})\in Q$.
Then
\begin{eqnarray*}
[y,l]u
&=&(yl)u-((y_{(-1)}\cdot \beta^{-1}(l))\beta(y_{0}))u\\
&=&(yl)u-\beta^{-1}(y_{(-1)}\cdot \beta^{-1}(l))(\beta(y_{0})\beta^{-1}(u))\\
&=&(yl)u-\beta^{-1}(y_{(-1)}\cdot \beta^{-1}(l))\beta^{-1}(\beta^{2}(y_{0})u)\\
&=&(yl)u=\beta(y)(l\beta^{-1}(u))\\
&=&\beta(y)[l,\beta^{-1}(u)]+\beta(y)((l_{(-1)}\cdot \beta^{-2}(u))\beta(l_{0}))\\
&=&\beta(y)[l,\beta^{-1}(u)]+(y(l_{(-1)}\cdot \beta^{-2}(u)))l_{0}\\
&=&\beta(y)[l,\beta^{-1}(u)].
\end{eqnarray*}
Since $\beta^{-1}(u)\in U,\beta(y)\in Q$, we obtain $[l,\beta^{-1}(u)]\in U$, $\beta(y)[l,\beta^{-1}(u)]=0$,
 and thus $[y,l]u$. Which means
$[Q,[L,L]]\subseteq Q$ and
$Q[L,L]\subseteq Q$. Hence
$$
L=QL=Q(Q+[Q,L])\subseteq Q.
$$
This implies $LU=0$, which contradicts the assumption $U\neq 0$. Hence, $Q=0$, and so
$ad^2_x(L)=0$. Similarly to case (1), one get $x\in Z_H(L)$.$\hfill \Box$
\medskip

\noindent{\bf Theorem 3.7.}
Let $(L, [,], \beta)$ be the derived braided Hom-Lie algebra.
Assume that char$(k)\neq 2$ and $L$ is $H$-simple.
If $V$ is an $H$-Hom-Lie ideal of $[L,L]$
such that any element in $V_0$ is $\beta$-invariant and
 $[V_0,V]\subseteq Z_H(L)_0$. Then $V_0\subseteq Z_H(L)_0$.
\medskip

{\bf Proof.}
For any $x\in V_0$. We consider the following two cases:

(1) If $ad_x(V)=0$, then $x\in Z_H(L)_0$ by Lemma 3.6.

(2) If $ad_x(V)\neq 0$, then for any $v\in V$ and $l\in L$, we have
\begin{eqnarray*}
[[x, [x, l]], v]
& \stackrel{(2.2)}{=} & -[[x,[x,l]]_{(-1)}\cdot\beta^{-1}(v),\beta([x,[x,l]]_0)]\\
& \stackrel{(2.1)}{=} & -[\beta(v_0),S^{-1}(v_{(-1)})\cdot\beta^{-1}([x,[x,l]])]\\
&=& -[\beta(v_0),\beta^{-1}(S^{-1}(\alpha(v_{(-1)}))\cdot[x,[x,l]])]\\
&=& -[\beta(v_0),\beta^{-1}([x,[x,S^{-1}(v_{(-1)})\cdot l]])]\\
&=& -[\beta(v_0),[x,[x,S^{-1}(\alpha^{-1}(v_{(-1)}))\cdot \beta^{-1}(l)]]].
\end{eqnarray*}
The fourth equality and the fifth equality hold since $x\in V_0$ is $\beta$-invariant.
By Lemma 3.3 (1), we get
\begin{eqnarray*}
&&(1\otimes C)(C\otimes 1)(v_0\otimes x\otimes[x,S^{-1}(\alpha^{-1}(v_{(-1)}))\cdot \beta^{-1}(l)])\\
&=&(1\otimes C)(x\otimes v_0\otimes [x,S^{-1}(\alpha^{-1}(v_{(-1)}))\cdot \beta^{-1}(l)])\\
&=&x\otimes v_{0(-1)}\cdot\beta^{-1}([x,S^{-1}(\alpha^{-1}(v_{(-1)}))\cdot \beta^{-1}(l)])\otimes \beta(v_{00})\\
&=&x\otimes v_{0(-1)}\cdot[x,S^{-1}(\alpha^{-2}(v_{(-1)}))\cdot \beta^{-2}(l)]\otimes \beta(v_{00})
   \end{eqnarray*}
   \begin{eqnarray*}
&=&x\otimes v_{(-1)2}\cdot[x,S^{-1}(\alpha^{-1}(v_{(-1)1}))\cdot \beta^{-2}(l)]\otimes v_{0}\\
&=&x\otimes [v_{(-1)21}\cdot x,v_{(-1)22}\cdot(S^{-1}(\alpha^{-1}(v_{(-1)1})\cdot \beta^{-2}(l))])\otimes v_{0}\\
&=&x\otimes [x,(\alpha^{-1}(v_{(-1)2})S^{-1}(\alpha^{-1}(v_{(-1)1})))\cdot \beta^{-2}(l)]\otimes v_{0}\\
&=&x\otimes [x,\epsilon(v_{(-1)})1\cdot \beta^{-2}(l)]\otimes v_{0}\\
&=&x\otimes [x,\beta^{-1}(l)]\otimes \beta^{-1}(v).
\end{eqnarray*}
Similarly,
$
(1\otimes C)(C\otimes 1)(v_0\otimes x\otimes[x,S^{-1}(\alpha^{-1}(v_{(-1)}))\cdot \beta^{-1}(l)])
= [x,\beta^{-1}(l)]\otimes \beta^{-1}(v)\otimes x.
$
By braided Hom-Jacobi identity, we have
\begin{eqnarray*}
[[x, [x, l]], v]
&=& -[\beta(v_0),[x,[x,S^{-1}(\alpha^{-1}(v_{(-1)}))\cdot \beta^{-1}(l)]]]\\
&=&[[\beta(x),l],[v,x]]+[\beta(x),[[x,\beta^{-1}(l)], \beta^{-1}(v)]]\\
&=&[[x,l],[v,x]]+[x,[[x,\beta^{-1}(l)], \beta^{-1}(v)]]\\
&\subseteq & [[x,L],[V,x]]+[x,[[x,L], \beta^{-1}(v)]]\\
&\subseteq & 0+[x,[[L,L],V]]\subseteq [x,V]\subseteq Z_H(L)_0.
\end{eqnarray*}
We obtain $[ad^2_x(L),V]\subseteq Z_H(L)_0$. By Lemma 3.5 (1), we have $ad^2_x(xl)=\beta^2(x)ad^2_x(l)$.

(2.1) If $ad^2_x(l)\neq 0$ for some $l\in L$, then $(ad^2_x(l))^{-1}\in Z_H(L)_0$ by Lemma 3.4. In this case,
it is easy to see that $x\in Z_H(L)_0$.

(2.2) Now we assume $ad^2_x(L)\varsubsetneq Z_H(L)_0$. Let $y\in L$
with $ad^2_x(y)\notin Z_H(L)_0$. Then we choose $z\in V$ such that $0 \neq ad_z(x)=u\in Z_H(L)_0$. Thus
there exist $v_1, v_2, v_3\in Z_H(L)_0$ such that $[z, ad^2_x(y)]=v_1$, $[\beta(z), ad^2_x(xy)]=v_2$ and
$[\beta^2(z), ad^2_x(x^2y)]=v_3$. Now we have
\begin{eqnarray*}
v_2
&=& [\beta(z), ad^2_x(xy)]=[\beta(z), xad^2_x(y)]\\
&=& [z,x]\beta(ad^2_x(y))+(\alpha(z_{(-1)})\cdot x)[\beta(z_{0}),ad^2_x(y)]\\
&=& [z,x]\beta(ad^2_x(y))+x[z,ad^2_x(y)]\\
&=& u\beta(ad^2_x(y))+xv_1.
\end{eqnarray*}
By Lemma 3.4, $u$ is invertible. Thus
$
ad^2_x(y)=\beta^{-1}(u^{-1}v_2-u^{-1}(xv_1)).
$
However, $v_1\in Z_H(L)$, $x\in V_0$, by Lemma 3.3 (1), we have
$
xv_1=v_1x,
$
and so
$
ad^2_x(y)=\beta^{-1}(u^{-1}v_2-u^{-1}(v_1x)).
$
Similarly, we have
\begin{eqnarray*}
v_3
&=& [\beta^2(z), ad^2_x(x^2y)]=[\beta(\beta(z)), xad^2_x(xy)]\\
&=& [\beta(z), x]\beta(ad^2_x(xy))+(\alpha((\beta(z))_{(-1)})\cdot x)[\beta((\beta(z))_{0}),ad^2_x(xy)]\\
&=& [\beta(z), x]\beta(ad^2_x(xy))+(\alpha^{2}(z_{(-1)})\cdot x)[\beta^{2}(z_{0}),ad^2_x(xy)]\\
&=& [\beta(z), \beta(x)]\beta(ad^2_x(xy))+x[\beta(z), ad^2_x(xy)]\\
&=& \beta(u)\beta(ad^2_x(xy))+xv_2\\
&=& u\beta(ad^2_x(xy))+xv_2.
\end{eqnarray*}
The last equality holds since  $u=ad_z(x)\in V_0$.
Thus
$
ad^2_x(xy)=\beta^{-1}(u^{-1}v_3-u^{-1}(v_2x)).
$
Using Lemma 3.5 (1), we have
\begin{eqnarray*}
ad^2_x(xy)
&=& xad^2_x(y)
=x\beta^{-1}(u^{-1}v_2-u^{-1}(v_1x))\\
&=& \beta^{-1}(\beta(x)(u^{-1}v_2)-\beta(x)(u^{-1}(v_1x)))\\
&=& \beta^{-1}((xu^{-1})\beta(v_2)-(xu^{-1})\beta(v_1x))\\
&=& \beta^{-1}((u^{-1}x)\beta(v_2)-(u^{-1}x)\beta(v_1x))\\
&=& \beta^{-1}(\beta(u^{-1})(xv_2)-\beta(u^{-1}x)(\beta(v_1)\beta(x)))\\
&=& \beta^{-1}(\beta(u^{-1})(v_2 x)-((u^{-1}x)\beta(v_1))\beta^{2}(x))\\
&=& \beta^{-1}((u^{-1}v_2)\beta(x)-(\beta(u^{-1})(xv_1))\beta^{2}(x))\\
&=& \beta^{-1}((u^{-1}v_2)\beta(x)-u^{-1}((xv_1)\beta(x)))\\
&=& \beta^{-1}(\beta(u^{-1})(v_2x)-u^{-1}((v_1x)\beta(x)))\\
&=& \beta^{-1}(u^{-1}(v_2x)-u^{-1}(\beta(v_1)x^{2})).
\end{eqnarray*}
Hence, $\beta(v_1)x^2-2v_2x+v_3=0$, that is, $x^2+\theta^1 x+ \theta^0=0$, where
$\theta^1=-2v_2/\beta(v_1)$, $\theta^0=v_3/\beta(v_1)$, and $\theta^1, \theta^0\in Z_H(L)$.
It is easy to see that $\theta^0=v_3/\beta(v_1)=(-\beta(v_1)x^2+2v_2x)/\beta(v_1)=-x^2-\theta^1x$.
 By Lemma 3.2 (2) and Lemma 3.3 (1) we have
\begin{eqnarray*}
0
&=&[-\theta^0,\beta(z)]=[x^2,\beta(z)]+[\theta^1 x,\beta(z)]\\
&=& \beta([x^2,z])+\beta(\theta^1)[x,z]+[\theta^1,x_{(-1)}\cdot\beta^{-1}(z)]\beta^{2}(x_0)\\
&=& \beta([x^2,z])+\beta(\theta^1)[x,z].
\end{eqnarray*}
By Lemma 3.3(1), one has $\beta([x^2,z])=-\beta(\theta^1)[x,z]=\beta(\theta^1)u.$
Similarly,
$$
\beta([x^2,z])=\beta(x[x,z]+[x,z]x)=2\beta([x,z]x)=-2\beta(ux)=-2ux.
$$
Since $u\in Z_H(L)_0$, $\beta(\theta^1)=-2x$, it follows that $\theta^1=-2\beta^{-1}(x)=-2x$.
 As char$(k)\neq 2$, we have $x=-(1/2)\theta^1\in Z_H(L)$, as desired.
$\hfill \Box$

\section{Universal enveloping  algebras of braided Hom-Lie\\ algebras}
\def\theequation{\arabic{section}.\arabic{equation}}
\setcounter{equation} {0}

In this section, we will first present the structure of the universal enveloping  algebra $U(L)$ of a braided Hom-Lie algebra $L$,
then we  show that $U(L)$ is a cocommutative  Hom-Hopf algebra.

\medskip

\medskip

\noindent{\bf Definition 4.1.}
 Let $(L, [,],\beta)$  be a braided Hom-Lie algebra.
 A universal enveloping algebra of $L$
  is a    monoidal Hom-algebra
   $$U(L)= (U(L),m_{U} ,\beta_{U} )$$
  together with a morphism $\psi: L\rightarrow  U(L)$ of Hom-Lie algebras in  $^{H}_{H}\mathcal{HYD}$  such that  the following
universal property holds: for any   monoidal Hom-algebra $A= (A,m_A ,\beta_A )$ and any Hom-Lie algebra
morphism $f: L\rightarrow A^{-}$ in  $^{H}_{H}\mathcal{HYD}$, there exists a unique morphism $g: U(L)\rightarrow A$ of  monoidal Hom-algebra in  $^{H}_{H}\mathcal{HYD}$ such that $g\circ \psi=f$.
\medskip

\noindent{\bf Definition 4.2.}
Let $(M,\beta_{M})$ be an involutive (i.e., $\beta_{M}^{2}=id$) Hom-Yetter-Drinfeld module.
A free involutive  monoidal Hom-algebra on $M$ is an involutive  monoidal Hom-algebra $(F_{M},\ast,\beta_{M})$
together with a morphism $j:M\rightarrow F_{M}$ in  $^{H}_{H}\mathcal{HYD}$,
satisfying  the following   property:   for any involutive  monoidal Hom-algebra $(A,\beta_{A})$ together with a morphism
 $f:M \rightarrow A$ in  $^{H}_{H}\mathcal{HYD}$,
  there is a unique morphism $\bar{f}:M\rightarrow F_{M}$ in  $^{H}_{H}\mathcal{HYD}$ such that $\bar{f}\circ j=f.$
 \medskip

 The well-known construction of the (non-unitary) free associative algebra on a module is
the tensor algebra equipped with the concatenation tensor product.
Recently, Guo, Zhang and Zheng generalized this method to Hom-associative algebras in \cite{Guo},
and  Armakan  Silvestov and Farhangdoost generalized the work to color Hom-associative algebras.
 Next we hope to extend the above work to  monoidal Hom-algebras  in  $^{H}_{H}\mathcal{HYD}$.
 \medskip

Let  $(M,\beta)$ be an  involutive  Hom-Yetter-Drinfeld module  and  $T(M)=\bigoplus_{i\geq 1}M^{\otimes i}$.
Obviously, $T(M)$  is an object in $^{H}_{H}\mathcal{HYD}$.
Define  the linear map $\beta_{T}$ and the binary operation $\odot$ on $T(M)$ as follows:
\begin{eqnarray*}
&&\beta_{T}(x)=\beta_{T}(x_{1}\otimes x_{2}\otimes \cdots \otimes x_{i})=\beta(x_{1})\otimes \beta(x_{2})\otimes \cdots \otimes \beta(x_{i}),\\
&&x\odot y=(x_{1}\otimes x_{2}\otimes \cdots \otimes x_{i})\odot(y_{1}\otimes y_{2}\otimes \cdots \otimes y_{j})
=\beta_{T}^{j-1}(x)\otimes y_{1}\otimes \beta_{T}(y_{2}\otimes \cdots \otimes y_{j}).
\end{eqnarray*}
One may check directly that $\beta_{T}$ and $\odot$ are morphisms in  $^{H}_{H}\mathcal{HYD}$.
Similar to the proof in \cite{Guo}, $(T(M),\odot,\beta_{T})$ is an involutive  monoidal Hom-algebra   in  $^{H}_{H}\mathcal{HYD}$.
 \medskip
%
%
%
%
%

\noindent{\bf Theorem 4.3.}
Let  $(H, \alpha)$ be an  involutive monoidal Hom-Hopf algebra and $(L, [,],\beta)$  an  involutive  braided Hom-Lie algebra.
Let $U(L)=T(L)/I$, where $I$ is the H-Hom-ideal of $T(L)$ generated by
$$\{x\o y-(x_{-1}\c \beta(y))\o \beta(x_{0})-[x,y]|~x,y\in L\}.$$
Let $\psi$ be the composition of the natural inclusion  $i:L\rightarrow T(L)$ with the canonical map $\pi: T(L)\rightarrow T(L)/I$.
Then $(U(L),\psi,\beta_{T})$ is an universal enveloping algebra of $L$.
 \medskip

{\bf Proof.}
We first show that $I$ is an object in $^{H}_{H}\mathcal{HYD}$. For any $x,y\in L$ and $ h\in H$, it is clear that
$
\rho(h_{1}\cdot x)
=(h_{111}\alpha^{-1}(x_{(-1)}))S(h_{12})\o\alpha(h_{112})\cdot x_{0}
=(\alpha^{-1}(h_{11})\alpha^{-1}(x_{(-1)}))S\alpha(h_{122})\o\alpha(h_{121})\cdot x_{0}.
$
Then we have
\begin{eqnarray*}
&& h\cdot (x\o y-(x_{-1}\c \beta(y))\o \beta(x_{0})-[x,y])\\
&=&h_1\cdot x\otimes h_2\cdot y- h_1\c (x_{-1}\c \beta(y))\otimes h_2\cdot\beta(x_{0}) -[h_1\cdot x, h_2\cdot y]  \\
&=&h_1\cdot x\otimes h_2\cdot y-(\alpha^{-1}(h_1)x_{-1})\c y\otimes h_2\cdot\beta(x_{0}) -[h_1\cdot x, h_2\cdot y]\\
&=&h_1\cdot x\otimes h_2\cdot y- (h_1\cdot x)_{-1}\c \beta(h_2\cdot y)\o \beta((h_1\cdot x)_{0})-[h_1\cdot x, h_2\cdot y]\in I.
\end{eqnarray*}
The last equality holds since
\begin{eqnarray*}
&& (h_1\cdot x)_{-1}\c \beta(h_2\cdot y)\o \beta((h_1\cdot x)_{0})\\
&=&((\alpha^{-1}(h_{11})\alpha^{-1}(x_{(-1)}))S\alpha(h_{122}))\c (\alpha(h_2)\cdot\beta(y))\o\alpha^{2}(h_{121})\cdot \beta(x_{0})\\
&=&(((\alpha^{-2}(h_{11})\alpha^{-2}(x_{(-1)}))S(h_{122}))\alpha(h_2))\cdot y\o\alpha^{2}(h_{121})\cdot \beta(x_{0})\\
&=&((\alpha^{-1}(h_{11})\alpha^{-1}(x_{(-1)}))(S(h_{122})h_2))\cdot y\o\alpha^{2}(h_{121})\cdot \beta(x_{0})\\
&=&((\alpha^{-2}(h_{1})\alpha^{-1}(x_{(-1)}))(S(h_{212})\alpha(h_{22})))\cdot y\o\alpha^{2}(h_{211})\cdot \beta(x_{0})\\
&=&((\alpha^{-2}(h_{1})\alpha^{-1}(x_{(-1)}))(S(h_{221})\alpha^{2}(h_{222})))\cdot y\o\alpha(h_{21})\cdot \beta(x_{0})\\
&=&((\alpha^{-2}(h_{1})\alpha^{-1}(x_{(-1)}))(\epsilon(h_{22})1_{H}))\cdot y\o\alpha(h_{21})\cdot \beta(x_{0})\\
&=&(\alpha^{-1}(h_{1})x_{(-1)})\cdot y\o h_{2}\cdot \beta(x_{0}).
\end{eqnarray*}
So $I$ is $H$-stable. Now we prove that $I$ is also $H$-costable, that is,
$\rho(x\o y-(x_{(-1)}\c \beta(y))\o \beta(x_{0})-[x,y])\in H\o I$, we note that
$
\rho(x_{(-1)}\c \beta(y))
=(x_{(-1)11}y_{(-1)})S(x_{(-1)2})\o\alpha(x_{(-1)12})\c \beta(y_0)
$
and compute
\begin{eqnarray*}
&&\rho(x_{-1}\c \beta(y)\o \beta(x_{0}))\\
&=&(x_{-1}\c \beta(y))_{(-1)}\alpha(x_{0(-1)})\o(x_{-1}\c \beta(y))_{0}\o\beta(x_{00})\\
&=&((x_{(-1)11}y_{(-1)})S(x_{(-1)2}))\alpha(x_{0(-1)})\o\alpha(x_{(-1)12})\c \beta(y_0)\o\beta(x_{00})\\
&=&((\alpha(x_{(-1)111})y_{(-1)})S(x_{(-1)12}))\alpha(x_{(-1)2})\o\alpha^{2}(x_{(-1)112})\c \beta(y_0)\o x_{0}\\
&=&((x_{(-1)11}y_{(-1)})S(x_{(-1)21}))\alpha^{2}(x_{(-1)22})\o\alpha(x_{(-1)12})\c \beta(y_0)\o\beta^{2}(x_{0})\\
&=&(\alpha(x_{(-1)11})\alpha(y_{(-1)}))(S(x_{(-1)21})\alpha(x_{(-1)22}))\o\alpha(x_{(-1)12})\c \beta(y_0)\o x_{0}\\
&=&(\alpha(x_{(-1)11})\alpha(y_{(-1)}))(\epsilon(x_{(-1)2})1_{H})\o\alpha(x_{(-1)12})\c \beta(y_0)\o x_{0}\\
&=&(\alpha^{2}(x_{(-1)1})\alpha(y_{(-1)}))1_{H}\o\alpha^{2}(x_{(-1)2})\c \beta(y_0)\o x_{0}\\
&=&\alpha(x_{(-1)1})y_{(-1)}\o x_{(-1)2}\c \beta(y_0)\o x_{0}\\
&=&x_{(-1)}y_{(-1)}\o x_{0(-1)}\c \beta(y_0)\o \beta(x_{00}).
\end{eqnarray*}
Therefore, we have
\begin{eqnarray*}
&&\rho(x\o y-(x_{(-1)}\c \beta(y))\o \beta(x_{0})-[x,y])\\
&=&x_{(-1)}y_{(-1)}\o x_{0}\o y_{0}-x_{(-1)}y_{(-1)}\o x_{0(-1)}\c \beta(y_0)\o \beta(x_{00})-x_{(-1)}y_{(-1)}\o [x_{0},y_{0}]\\
&=&x_{(-1)}y_{(-1)}\o (x_{0}\o y_{0}-x_{0(-1)}\c \beta(y_0)\o \beta(x_{00})-[x_{0},y_{0}])\in H\o I,
\end{eqnarray*}
as desired, where $\rho[x,y]=x_{(-1)}y_{(-1)}\o [x_{0},y_{0}]$ since $[, ]$ is a morphism in  $^{H}_{H}\mathcal{HYD}$.

Next, we show that $\psi$ is a morphism of braided Hom-Lie algebras.
It is easy to see that $\psi$ is a morphism in  $^{H}_{H}\mathcal{HYD}$.
Now we prove that $\psi$ is compatible with the bracket product, we denote the multiplication in $U(L)$ by $\ast$ and calculate
\begin{eqnarray*}
\psi([x,y])
&=&\pi([x,y])=\pi(x\o y-(x_{(-1)}\c \beta(y))\o \beta(x_{0}))\\
&=&\pi(x\odot y-(x_{(-1)}\c \beta(y))\odot \beta(x_{0}))\\
&=&\pi(x)\ast\pi(y)-\pi(x_{(-1)}\c \beta(y))\ast\pi(\beta(x_{0}))\\
&=&\psi(x)\ast\psi(y)-\psi(x_{(-1)}\c \beta(y))\ast\psi(\beta(x_{0}))\\
&=&\psi(x)\ast\psi(y)-(x_{(-1)}\c \psi(\beta(y)))\ast\psi(\beta(x_{0}))\\
&=&\psi(x)\ast\psi(y)-((\psi(x))_{(-1)}\c \beta(\psi(y)))\ast\beta((\psi(x))_{0})\\
&=&[\psi(x),\psi(y)].
\end{eqnarray*}

Finally, we show that the following statement holds:
for any involutive monoidal Hom-algebra of $(A,m_{A},\beta_{A})$
and any  homomorphism $f: L\longrightarrow A^{-1}$ of   Hom-Lie algebras in $^{H}_{H}\mathcal{HYD}$,
 there exists a unique morphism $g: U(L)\longrightarrow A$ in $^{H}_{H}\mathcal{HYD}$
 such that the following diagram commutes:
 $$
 \begin{array}{ccc}
 L&\stackrel {\psi}{\longrightarrow }&U(L) \\
 f\downarrow & \swarrow g& \, \\
 A &\, &
\end{array}
 $$
 To prove this statement, we first consider a unique homomorphism $f^*$ of $T(L)$ which maps
 $T(L)$ into $A$ by extending the homomorphism $f$ of $L$ into $A$. For any $x,y\in L$, we have
\begin{eqnarray*}
&&f^*(x\o y-(x_{(-1)}\c \beta(y))\o \beta(x_{0}))\\
&=&f^*(x\odot y-(x_{(-1)}\c \beta(y))\odot \beta(x_{0}))\\
&=&f^*(x)f^*(y)-f^*(x_{(-1)}\c \beta(y))f^*(\beta(x_{0}))\\
&=&f(x)f(y)-f(x_{(-1)}\c \beta(y))f(\beta(x_{0}))\\
&=&f(x)f(y)-x_{(-1)}\c \beta(f(y))\beta(f(x_{0}))\\
&=&[f(x),f(y)]=f([x,y])=f^*([x,y]).
\end{eqnarray*}
 This shows that $I\subset ker f^*$, and we have a unique homomorphism $g$ of $U(L)=T(L)/I$ into $A$ such that
 $g(x+I)=f(x)$  or  $g\psi(x)=f(x)$.
 Hence $f=g\psi$, since $L$ generates $T(L)$.

Furthermore, it is easy to see  that $\alpha_{A}\circ g=g\circ\beta_{T}$.
We still need to check that  $g$ is a morphism in $^{H}_{H}\mathcal{HYD}$.
Since  $\rho _Af=(1\otimes f)\rho _L$ by our assumption,
where $\rho _A$ and $\rho _L$ are the $(H, \alpha)$-Hom-comodule structure of $A$ and $L$ respectively,
 for any $\overline {x},\overline {y}\in U(L)$, we have
\begin{eqnarray*}
\rho _Ag(\overline {x}\ast\overline {y})
&=&\rho _A(g(\overline x)g(\overline y))=\rho _A(f(x)f(y))\\
&=&(f(x))_{(-1)}(f(y))_{(-1)}\otimes (f(x))_{0}(f(x))_{0}\\
&=&x_{(-1)}y_{(-1)}\otimes f(x_0)f(y_0)
=x_{(-1)}y_{(-1)}\otimes g(\overline{x_0})f(\overline{y_0})\\
&=&(1\otimes g) (x_{(-1)}y_{(-1)}\otimes(\overline{x_{0}}\ast\overline{y_{0}}))
=(1\otimes g)\rho _U(\overline{x}\ast\overline{y}),
\end{eqnarray*}
It follows that $g$ is indeed  $(H, \alpha)$-linear.
Similarly, one may check that $g$ is also $(H, \alpha)$-linear.
 And the proof is completed.
 \hfill $\square$

Now we will define a  Hom-Hopf algebra structure on the  universal enveloping  algebra $U(L)$,
we first present a useful Lemma.
\medskip

{\bf Lemma 4.4.}
Let  $(H, \alpha)$ be an  involutive monoidal  Hom-Hopf algebra and $(L, [,],\beta)$  an  involutive  braided Hom-Lie algebra.
Assume $U(L)$ is the universal enveloping   algebra of $L$.
Then there exists a homomorphism
$g: U(L\oplus L)\longrightarrow U(L)\otimes U(L)$
of monoidal Hom-algebras in $^{H}_{H}\mathcal{HYD}$.
\medskip

{\bf Proof.} Define $f: L \oplus L\longrightarrow U(L)\otimes U(L)$
 by
 $$(x, y)\mapsto \beta_{T}(\overline x) \otimes 1+1\otimes\beta_{T}(\overline y).$$

 We first show that $f$ is a  morphism in $^{H}_{H}\mathcal{HYD}$.
 In fact, for any $h\in H$ and $x,y\in L$, we have
\begin{eqnarray*}
h\cdot f(x, y))
   &=&h_{1}\cdot \beta_{T}(\overline x)\otimes h_{2}\cdot 1+h_{1}\cdot 1\otimes h_{2}\cdot\beta_{T}(\overline y)\\
   &=&h_{1}\cdot \beta_{T}(\overline x)\otimes \epsilon(h_{2})1+\epsilon(h_{1})1\otimes h_{2}\cdot\beta_{T}(\overline y)\\
   &=&\alpha(h)\cdot \beta_{T}(\overline x)\otimes 1+ 1\otimes \alpha(h)\cdot\beta_{T}(\overline y)\\
   &=&\beta_{T}(h\cdot\overline x)\otimes 1+ 1\otimes \beta_{T}(h\cdot\overline y)\\
   &=&\beta_{T}(\overline{h\cdot x})\otimes 1+ 1\otimes  \beta_{T}(\overline{h\cdot y})\\
    &=&f(h\cdot x,h\cdot y)
    =f(h\cdot (x, y)).
\end{eqnarray*}
 It follows that $f$ is $H$-linear.
 Similarly, one may check that $f$ is $H$-colinear.

 Second, we prove that $f$ is a Hom-Lie homomorphism.
 For any $x,y',x,y'\in L$, we have
 \begin{eqnarray*}
[f(x,y),f(x',y')]]
&=&[\beta_{T}(\overline x) \otimes 1+1\otimes\beta_{T}(\overline y), \beta_{T}(\overline {x'}) \otimes 1+1\otimes\beta_{T}(\overline {y'})]\\
&=&[\beta_{T}(\overline x) \otimes 1,\beta_{T}(\overline {x'}) \otimes 1]
  +[\beta_{T}(\overline x) \otimes 1,1\otimes\beta_{T}(\overline{y'})]+\\
 && [1\otimes\beta_{T}(\overline y),\beta_{T}(\overline {x'}) \otimes 1]
   +[1\otimes\beta_{T}(\overline y),1\otimes\beta_{T}(\overline {y'})].
\end{eqnarray*}
Recall that multiplication in $U(L)\otimes U(L)$ is
 \begin{eqnarray*}
(\overline x\o\overline y)(\overline {x'}\o\overline {y'})
=\overline x(y_{(-1)}\c\beta_{T}^{-1}(\overline {x'}))\o(\beta_{T}(y_{0})y').
\end{eqnarray*}
Obviously, we have
$(\overline x\o 1)(1\o\overline y)=\beta_{T}(\overline x) \otimes\beta_{T}(\overline y)$  and
$(1\o\overline x)(\overline y\o 1)=\alpha(x_{(-1)})\c\overline y\o x_{0}.$
Therefore,
 \begin{eqnarray*}
[\beta_{T}(\overline x) \otimes 1,1\otimes\beta_{T}(\overline {y'})]
&=&(\beta_{T}(\overline x) \otimes 1)(1\otimes\beta_{T}(\overline {y'}))
   -((\alpha(x_{(-1)})1)\c(1\otimes\overline {y'}))(\overline x_{0}\otimes 1)\\
&=&\overline x\otimes\overline {y'}
   -(x_{(-1)}\c(1\otimes\overline {y'}))(\overline x_{0}\otimes 1)\\
&=&\overline x\otimes\overline {y'}
   -(1\otimes\alpha(x_{(-1)})\c\overline {y'}))(\overline x_{0}\otimes 1)\\
&=&\overline x\otimes\overline {y'}
   -((\alpha^{2}(x_{(-1)11})y_{(-1)})S\alpha(x_{(-1)2}))\c  \overline {x_{0}}\o x_{(-1)12}\c  \overline {y_{0}}\\
&=&\overline x\otimes\overline {y'}-\overline x\otimes\overline {y'}
=0,
\end{eqnarray*}
where $((\alpha^{2}(x_{(-1)11})y_{(-1)})S\alpha(x_{(-1)2}))\c  \overline {x_{0}}\o x_{(-1)12}\c  \overline {y_{0}}
=\overline x\otimes\overline {y'}$ since  the braiding is symmetric on $L$.
 Similarly, we have $[1\otimes\beta_{T}(\overline y),1\otimes\beta_{T}(\overline {y'})]=0.$
 Also,
 \begin{eqnarray*}
[\beta_{T}(\overline x) \otimes 1,\beta_{T}(\overline {x'}) \otimes 1]
&=&(\beta_{T}(\overline x)(1\c\overline {x'}))\otimes \beta_{T}(1)1-((\alpha(x_{(-1)})1)\c(\overline {x'}\o 1))(\overline {x_{0}}\o 1)\\
&=& \beta_{T}(\overline x)\beta_{T}(\overline {x'})\otimes 1-(\alpha(x_{(-1)})\c \overline {x'}\o 1)(\overline {x_{0}}\o 1)\\
&=& \beta_{T}(\overline x)\beta_{T}(\overline {x'})\otimes 1-(\alpha(x_{(-1)})\c \overline {x'})\overline {x_{0}}\o 1\\
&=& \beta_{T}(\overline x)\beta_{T}(\overline {x'})\otimes 1
    -((\beta_{T}(\overline {x}))_{(-1)}\c \beta_{T}^{-1}(\beta_{T}(\overline {x'})))\beta_{T}((\beta_{T}(\overline {x}))_{0})\o 1\\
&=&[\beta_{T}(\overline x),\beta_{T}(\overline {x'})]\o 1.
\end{eqnarray*}
 Similarly, we have
 $[1\otimes\beta_{T}(\overline y),1\otimes\beta_{T}(\overline {y'})]
 =1\o[\beta_{T}(\overline y),\beta_{T}(\overline {y'})].$
 Then we have
 \begin{eqnarray*}
[f(x,y),f(x',y')]]
&=&[\beta_{T}(\overline x),\beta_{T}(\overline {x'})]\o 1+1\o[\beta_{T}(\overline y),\beta_{T}(\overline {y'})]\\
&=&\beta_{T}([\overline x,\overline {x'}])\o 1+1\o\beta_{T}([\overline y,\overline {y'}])\\
&=&f([(x,y),(x',y')]).
\end{eqnarray*}
So $f$ is a Hom-Lie homomorphism.
 Now by the universal property of $U(L\oplus L)$, there exists a homomorphism
$g: U(L\oplus L)\longrightarrow U(L)\otimes U(L)$ of  monoidal Hom-algebras in $^{H}_{H}\mathcal{HYD}$.
\medskip

 {\bf Theorem 4.5.} Let  $(H, \alpha)$ be an  involutive monoidal Hom-Hopf algebra
  and $(L, [,],\beta)$  an  involutive  braided Hom-Lie algebra.
Then $U(L)$  in Theorem 4.3 is a monoidal Hom-Hopf algebra in $^{H}_{H}\mathcal{HYD}$ with
\begin{eqnarray*}
 &&\Delta (\overline l)=\beta_{T}(\overline l) \otimes 1+1\otimes\beta_{T}(\overline l);\\
 &&\Delta (1)=1\otimes 1,~~\epsilon (\overline l)=0, ~~\epsilon (1)=1;\\
&&S(\overline l)=-\overline l,~~S(\bar x  \bar y) =( x_{(-1)}\cdot S(\beta_{T}^{-1}(\overline y)))S(\beta_{T}(\overline {x_0})).
\end{eqnarray*}
for all $l\in L$ and  $\overline x, \overline y\in U(L)$.
\medskip

{\bf Proof.}    We first consider the diagonal
mapping $d: L\longrightarrow L\oplus L$  defined by $l\mapsto (l,l)$.
 It is easy to check that $d$ is a Hom-Lie homomorphism in  $^{H}_{H}\mathcal{HYD}$.
Let $f$ be the map described in Lemma 4.4.
Then $f\circ d$  is  a Hom-Lie homomorphism from $L$ to $ U(L)\o U(L)$,
therefore there exists a  homomorphism
$\Delta : U(L)\rightarrow U(L)\o U(L)$,
which is  a homomorphism of  monoidal Hom-algebras in  $^{H}_{H}\mathcal{HYD}$ satisfying the following condition
 $$
 \Delta (\overline l)=((f\circ d)(l))=\beta_{T}(\overline l) \otimes 1+1\otimes\beta_{T}(\overline l),
 $$
for all $\bar l\in \bar L.$
It is now straightforward to check that
$
(\beta_{T}^{-1}\otimes \Delta )\Delta =(\Delta \otimes \beta_{T}^{-1})\Delta
$
 and
 $
 (\e \otimes \beta_{T})\Delta =(\beta_{T}\otimes \epsilon ) \Delta =\beta_{T}^{-1}.
$

It is easy to see that $S$  is a well-defined  morphism in $^{H}_{H}\mathcal{HYD}$,
since if we define $\widetilde{S}$  on the free generators of $T(L)$ by
$\widetilde{S}(\overline l)=-\overline l, \widetilde{S}(1)=1,$
and set $\widetilde{S}(\bar x  \bar y) =( x_{(-1)}\cdot \widetilde{S}(\beta_{T}^{-1}(\overline y)))\widetilde{S}(\beta_{T}(\overline {x_0})),$
 then $\widetilde{S}$  is a morphism in $^{H}_{H}\mathcal{HYD}$ which vanishes on $I$. Thus $S$ is well defined.

To show that $S$  is an antipode, we first note that
\begin{eqnarray*}
(m(id\o S)\circ\Delta)(\overline l)
&=&m(id\o S)(\beta_{T}(\overline l) \otimes 1+1\otimes\beta_{T}(\overline l))\\
&=&m(\beta_{T}(\overline l) \otimes 1-1\otimes\beta_{T}(\overline l))
=0=\epsilon (\overline l),\\
(m(S\o id )\circ\Delta)(\overline l)
&=&m(S\o id)(\beta_{T}(\overline l) \otimes 1+1\otimes\beta_{T}(\overline l))\\
&=&m(-\beta_{T}(\overline l) \otimes 1+1\otimes\beta_{T}(\overline l))=0=\epsilon (\overline l),
\end{eqnarray*}
 for any generator  $l\in L$.
Similarly, one may check that $(m(id\o S)\circ\Delta)(1)=(m(S\o id )\circ\Delta)(1)=\epsilon (1)$.
Therefore,  we can derive that
\begin{eqnarray*}
(m(id\o S)\circ\Delta)(\overline{x} ~\overline{y})
&=&m(id\o S)(\overline{x_{1}}(x_{2(-1)}\c\beta_{T}^{-1}(\overline {y_{1}}))\o\beta_{T}(\overline {x_{20}})\overline {y_{2}})\\
&=&m(\overline{x_{1}}(x_{2(-1)}\c\beta_{T}^{-1}(\overline {y_{1}}))\o S(\beta_{T}(\overline {x_{20}})\overline {y_{2}}))\\
&=&\{\overline{x_{1}}(x_{2(-1)}\c\beta_{T}^{-1}(\overline {y_{1}}))\}\{(\alpha(x_{20(-1)})\c S\beta_{T}(\overline {y_{2}}))S(\overline {x_{200}})\}\\
&=&\{(\overline{x_{1}}(\alpha(x_{2(-1)1})\c\beta_{T}(\overline {y_{1}}))\}\{(\alpha(x_{2(-1)2})\c S\beta_{T}(\overline {y_{2}}))S\beta_{T}(\overline {x_{20}})\}\\
&=&\{\overline{x_{1}}\beta_{T}(x_{2(-1)1}\c\overline {y_{1}})\}\beta_{T}((x_{2(-1)2}\c S(\overline {y_{2}}))S(\overline {x_{20}}))\\
&=&\beta_{T}(\overline{x_{1}})(\beta_{T}(x_{2(-1)1}\c\overline {y_{1}})\{x_{2(-1)2}\c S(\overline {y_{2}}))S(\overline {x_{20}})\})\\
&=&\beta_{T}(\overline{x_{1}})(\{(x_{2(-1)1}\c\overline {y_{1}})(x_{2(-1)2}\c S(\overline {y_{2}}))\}S\beta_{T}(\overline {x_{20}}))\\
&=&\beta_{T}(\overline{x_{1}})\{(x_{2(-1)}\c\epsilon(\overline{y})1)S\beta_{T}(\overline {x_{20}})\}\\
&=&\epsilon(\overline{y})\beta_{T}(\overline{x_{1}})S\beta_{T}(\overline {x_{2}})
=\epsilon(\overline{y})\epsilon(\overline{x}).
\end{eqnarray*}
Similarly, we can show that $(m(S\o id)\circ\Delta)(\overline{x} ~\overline{y})=\epsilon(\overline{y})\epsilon(\overline{x})$.
So $S$ is an antipode on $U(L)$, and this finishes the proof.
 $\hfill \Box$
\medskip

 {\bf Corollary 4.6.}
Under the hypotheses of the Theorem 4.5,  the universal enveloping    algebra $U(L)$ is $H$-cocommutative.
\medskip

{\bf Proof.}
For any $\overline{x}\in U(L)$, we have
$
C_{U,U}\Delta(\overline{x} )
=C_{U,U}(\beta_{T}(\overline x) \otimes 1+1\otimes\beta_{T}(\overline x))
=\alpha(x_{(-1)})\c\beta_{T}^{-1}(1)\otimes\beta_{T}^{2}(\overline x_{0})+1\c\beta_{T}^{-1}\beta_{T}(\overline x)\otimes\beta_{T}(1)
=1\otimes\beta_{T}(\overline x))+\beta_{T}(\overline x) \otimes 1
=\Delta(\overline{x} ).
$
It follows that
$
C_{U,U}\Delta=\Delta,
$
as desired.
 $\hfill \Box$
 \medskip

As an application of Theorem 4.5,   we will define a Hom-Yetter-Drinfeld module structure on the $End(V)$ and construct a Radford's Hom-biproduct. In order to define a good $(H,\a)$-Hom-module operation on $End(V)$,
it is necessary to assume that $\a=id_{H}.$

  {\bf Lemma 4.7.}
 Let $H$ be a Hopf algebra with a bijective antipode and $(V,\nu)$  a finite-dimensional Hom-Yetter-Drinfeld module in  $^{H}_{H}\mathcal{HYD}$.
 Then  $(End(V),\delta)$ is a Hom-Yetter-Drinfeld module under the following structures
 \begin{eqnarray*}
&&(h\c f)(v)=h_1\c f(S(h_2)\c v),~\delta(f)(v)=f(\nu^{2}(v)),\\
&& \rho(f)(v)=(f(v_0))_{(-1)}S^{-1}(v_{(-1)})\o(f(v_0))_{0},
\end{eqnarray*}
for any $v\in V$.
 \medskip

 {\bf Proof.}
 We first show that $(End(V),\delta)$ is a Hom-module.
 In fact, for any $h,g\in H$, $f\in End(V)$ and $v\in V$, we have
  \begin{eqnarray*}
(h\c(g\c f))(v)
&=&h_1\c (g\c f)(S(h_2)\c v)
=h_1\c (g_1\c f(S(g_2)\c (S(h_2)\c v)))\\
&=&h_1\c (g_1\c f(S(g_2)S(h_2)\c \nu(v)))
=(h_1g_1)\c f(S(g_2)S(h_2)\c \nu^{2}(v)),\\
((hg)\c\delta(f))(v)
&=&(hg)_{1}\c\delta(f)(S((hg)_2))\c v)
=(h_1g_1)\c f(S(h_2g_2)\c \nu^{2}(v)).
\end{eqnarray*}
 It follows that $h\c(g\c f)=(hg)\c\delta(f)$.
Now we verify $1_H \cdot f=\delta(f)$ and $\delta(h\cdot f)=h\cdot\delta(f)$ as follows
\begin{eqnarray*}
(1_H \cdot f)(v)
&=&1\c f(1\c v)
=1\c f(\nu(v))
=f(\nu^{2}(v))\\
\delta(h\cdot f)(v)
&=&(h\cdot f)(\nu^{2}(v))
=h_1\c f(S(h_2)\c \nu^{2}(v))\\
&=&h_1\c \delta(f)(S(h_2)\c v)
=(h\c\delta(f))(v).
\end{eqnarray*}
 So $(End(V),\delta)$ is a Hom-module, as desired.
 Similarly, one may check that $(End(V),\delta)$ is a Hom-comodule.

Now we show that for any $f\in End(V) $ and  $h\in H$, the
 following compatibility condition
$$
h_1 f_{(-1)}\otimes h_2\cdot f_0=(h_1\cdot \delta^{-1}(f))_{(-1)}h_2
\otimes \delta((h_1\cdot \delta^{-1}(f))_{0}),
$$
holds.
For this,  we take  $h\in H, f\in End(V) , v\in V$. On the one hand, we have
\begin{eqnarray*}
 &&(h_1\cdot \delta^{-1}(f))_{(-1)}h_2\otimes \delta((h_1\cdot \delta^{-1}(f))_{0})(v) \\
 &=&(h_1\cdot \delta^{-1}(f))_{(-1)}h_2\otimes(h_1\cdot \delta^{-1}(f))_{0}(\nu^{2}(v))\\
 &=& ((h_1\cdot \delta^{-1}(f))(\nu^{2}(v_{00})))_{(-1)}S^{-1}(v_{(-1)})h_2\o((h_1\cdot \delta^{-1}(f))(\nu^{2}(v_{00})))_{0}\\
  &=& (h_1\cdot f(S(h_3)\cdot v_0))_{(-1)}S^{-1}(v_{(-1)})h_3\o(h_1\cdot f(S(h_3)\cdot v_0))_{0}\\
    &=&h_1 (f(S(h_4)\cdot v_0))_{(-1)}S(h_3)S^{-1}(v_{(-1)})h_5\o h_3\cdot(f(S(h_4)\cdot v_0))_{0}.
\end{eqnarray*}
On the other hand, we have
\begin{eqnarray*}
 & & h_1f_{(-1)}\otimes (h_2\cdot f_0)(v) \\
  &=&h_1 f_{(-1)}\otimes h_2\cdot (f_0((S(h_3))\cdot v))  \\
  &=&h_1(f(S(h_3))\cdot v)_0)_{(-1)}S^{-1}(S(h_3)\cdot v)_{(-1)}
  \otimes h_2\cdot (f(((Sh_3))\cdot v)_0)_0   \\
  &=&h_1(f(S(h_4)\cdot v_0))_{(-1)}S^{-1}(S(h_5)v_{(-1)}S^{2}h_3)
  \otimes h_2\cdot (f(S(h_4))\cdot v_0))_0 \\
  &=&h_1(f(S(h_4)\cdot v_0))_{(-1)}S(h_3)S^{-1}(v_{(-1)})h_5
  \otimes h_2\cdot (f((Sh_4))\cdot v_0)_0.
\end{eqnarray*}
 So $(End(V),\delta)\in ^{H}_{H}\mathcal{HYD}$. The proof is finished.
 $\hfill \Box$
 \medskip

   {\bf Lemma 4.8.}
 Let $H$ be a Hopf algebra with a bijective antipode and $(V,\nu)$  a finite-dimensional involutive  Hom-Yetter-Drinfeld module in  $^{H}_{H}\mathcal{HYD}$.
 Then  $(End(V),\delta)$ is a monoidal Hom-algerba in $^{H}_{H}\mathcal{HYD}$.
 \medskip

 {\bf Proof.}
We first show that  $(V,\nu)$  is a $H$-module algerba.
Indeed, for any  $h\in H, f,g\in End(V) $ and $ v\in V$,  we have
\begin{eqnarray*}
(h_1\c f)(h_2\c g)(v)
&=&(h_1\c f)(h_2\c g(S(h_3)\c v))\\
&=&h_1\c f(S(h_2)\c(h_3\c g(S(h_4)\c v)))\\
&=&h_1\c f((S(h_2)h_3)\c g(S(h_4)\c \nu(v)))\\
&=&h_1\c f((\epsilon(h_2)1_H)\c g(S(h_3)\c \nu(v)))\\
&=&h_1\c f(g(S(h_2)\c \nu^2(v)))\\
&=&h_1\c (fg)(S(h_2)\c v).
\end{eqnarray*}
It follows that $h\c(fg)=(h_1\c f)(h_2\c g)$. Also, we have
 \begin{eqnarray*}
(h\c id)(v)
&=&h_1\c id(S(h_2)\c v)
=h_1\c(S(h_2)\c v)\\
&=&(h_1S(h_2))\c \nu(v)
=\epsilon(h)1_H\c \nu(v)
=\epsilon(h)v.
\end{eqnarray*}
 So $ h\c id=\epsilon(h)id$. Therefore,  $(V,\nu)$  is a $H$-module algerba.

 Next, we will show that $(V,\nu)$  is a $H$-comodule algerba.
In fact, for any   $ f,g\in End(V)$ and $ v\in V$,  we have
  \begin{eqnarray*}
(fg)_{(-1)}\o(fg)_{0}(v)
&=&((fg)(v_0))_{(-1)}S^{-1}(v_{(-1)})\o((fg)(v_0))_{0}\\
&=&(fg(v_0))_{(-1)}S^{-1}(v_{(-1)})\o(fg(v_0))_{0},\\
f_{(-1)}g_{(-1)}\o f_{0}g_{0}(v)
&=&f_{(-1)}(g(v_0))_{(-1)}S^{-1}(v_{(-1)})\o f_{0}((g(v_0))_{0})\\
&=&(f((g(v_0))_{00}))_{(-1)}S^{-1}((g(v_0))_{0(-1)})(g(v_0))_{(-1)}S^{-1}(v_{(-1)})\\
  &&\o(f((g(v_0))_{00}))_{0}\\
&=&(f(\nu^{-1}(g(v_0))_{0}))_{(-1)}S^{-1}((g(v_0))_{(-1)2})(g(v_0))_{(-1)1}S^{-1}(v_{(-1)})\\
  &&\o(f(\nu^{-1}(g(v_0))_{0}))_{0}\\
&=&(f(\nu^{-1}(g(v_0))_{0}))_{(-1)}\epsilon(g(v_0)_{(-1)})S^{-1}(v_{(-1)})\o(f(\nu^{-1}(g(v_0))_{0}))_{0}\\
&=&(f((g(v_0))_{0}))_{(-1)}S^{-1}(v_{(-1)})\o(f((g(v_0))_{0}))_{0}\\
&=&(fg(v_0))_{(-1)}S^{-1}(v_{(-1)})\o(fg(v_0))_{0}.
\end{eqnarray*}
It follows that $(fg)_{(-1)}\o(fg)_{0}=f_{(-1)}g_{(-1)}\o f_{0}g_{0}.$
 Also, we have
   \begin{eqnarray*}
 \rho(id)(v)
 &=&v_{0(-1)}S^{-1}(v_{(-1)})\o v_{00}
 =v_{(-1)2}S^{-1}(v_{(-1)1})\o \nu^{-1}(v_{0})\\
 &=&\epsilon(v_{(-1)})1_{H}\o \nu^{-1}(v_{0})
 =1_{H}\o v=1_{H}\o id(v).
\end{eqnarray*}
So  $\rho(id)=1_{H}\o id$, as desired. And this complete the proof.
 $\hfill \Box$
 \medskip

   {\bf Lemma 4.9.}
 Let $H$ be a Hopf algebra with a bijective antipode and $(V,\nu)$  a finite-dimensional involutive  Hom-Yetter-Drinfeld module in  $^{H}_{H}\mathcal{HYD}$. Assume that the braiding $C$ is symmetric on $V$.
 Then  $(End(V),\delta)$ is a braided Hom-Lie algebra, where the bracket product is defined by
 $$[f,g]=fg-(f_{(-1)}\c\delta^{-1}(g))\delta(f_0),$$
 for any $f,g\in End(V).$
 \medskip

 {\bf Proof.}
 Since the braiding $C$ is symmetric on $V$, one may check that   $C$ is symmetric on $End(V)$, too.
By Proposition 2.2, $(End(V),\delta)$ is a braided Hom-Lie algebra.
 $\hfill \Box$
 \medskip

{\bf Proposition 4.10.}
Let $H$ be a Hopf algebra with a bijective antipode and $(V,\nu)$  a finite-dimensional involutive  Hom-Yetter-Drinfeld module.
 Assume that the braiding $C$ is symmetric on $V$.
 Then the Radford's Hom-biproduct  $(U(End(V))_{\sharp}^{\times} H, \delta\o id)$ is a monoidal Hom-Hopf algebra,
 where the multiplication is defined by
$$(f\times h)(f'\times h')=f(h_1\c \delta^{-1}(f))\times h_2h',$$
 the coproduct  is defined by
 $$\Delta (f\times h)=(f_1\times f_{2(-1)}h_1)\otimes(\delta(f_{2(0)})\times h_2),$$
 the antipode is defined by
$$S(f\times h)=(1\times S(f_{(-1)}h))(S(f_0)\times 1),$$
for all $f\times h, f'\times h'\in U(End(V))_{\sharp}^{\times}  H.$
 \medskip

 {\bf Proof.}
 By  Lemma 4.9 and Theorem 4.5, $(U(End(V)),\delta)$ is a  monoidal  Hom-Hopf algebra in $^{H}_{H}\mathcal{HYD}$.
 By Proposition 3.6 in \cite{LIU2014}, $(U(End(V))_{\sharp}^{\times} H, \delta\o id)$ is a  monoidal  Hom-Hopf algebra.

\begin{center}
 {\bf ACKNOWLEDGEMENT}
 \end{center}

   The work of S. X. Wang is  supported by  the outstanding top-notch talent cultivation project of Anhui Province (No. gxfx2017123)
 and the Anhui Provincial Natural Science Foundation (1808085MA14).
 The work of X. H. Zhang is  supported by the NSF of China (No. 11801304, 11801306) and the Project Funded by China Postdoctoral Science Foundation (No. 2018M630768).
The work of S. J. Guo is  supported by the NSF of China (No. 11761017) and the Youth Project for Natural Science Foundation of Guizhou provincial department of education (No. KY[2018]155).

\renewcommand{\refname}{REFERENCES}

\end{document}